# QUASI-MAXIMUM-LIKELIHOOD ESTIMATION IN CONDITIONALLY HETEROSCEDASTIC TIME SERIES: A STOCHASTIC RECURRENCE EQUATIONS APPROACH[1]

By Daniel Straumann and Thomas Mikosch

*RiskMetrics Group and University of Copenhagen*

This paper studies the quasi-maximum-likelihood estimator (QMLE) in a general conditionally heteroscedastic time series model of multiplicative form $X_t = \sigma_t Z_t$, where the unobservable volatility $\sigma_t$ is a parametric function of $(X_{t-1}, \ldots, X_{t-p}, \sigma_{t-1}, \ldots, \sigma_{t-q})$ for some $p, q \geq 0$, and $(Z_t)$ is standardized i.i.d. noise. We assume that these models are solutions to stochastic recurrence equations which satisfy a contraction (random Lipschitz coefficient) property. These assumptions are satisfied for the popular GARCH, asymmetric GARCH and exponential GARCH processes. Exploiting the contraction property, we give conditions for the existence and uniqueness of a strictly stationary solution $(X_t)$ to the stochastic recurrence equation and establish consistency and asymptotic normality of the QMLE. We also discuss the problem of invertibility of such time series models.

**1. Introduction.** Gaussian quasi-maximum-likelihood estimation, that is, likelihood estimation under the hypothesis of Gaussian innovations, is a popular method which is widely used for inference in time series models. However, it is often a nontrivial task to establish the consistency and asymptotic normality of the quasi-maximum-likelihood estimator (QMLE) applied to specific models and, therefore, an in-depth analysis of the probabilistic structure generated by the model is called for. A classical example of this kind is the seminal paper by Hannan [18] on estimation in linear ARMA time series.

In this paper we study the QMLE for a general class of conditionally heteroscedastic time series models, which includes GARCH, asymmetric

Received January 2003; revised January 2006.

[1]Supported by Danish Research Council (SNF) Grant 21-01-0546.

GARCH and exponential GARCH. Recall that a GARCH$(p, q)$ [generalized autoregressive conditionally heteroscedastic of order $(p, q)$] process [4] is defined by

$$(1.1) \qquad X_t = \sigma_t Z_t, \qquad t \in \mathbb{Z},$$

and the recurrence equation

$$(1.2) \qquad \sigma_t^2 = \alpha_0 + \sum_{i=1}^{p} \alpha_i X_{t-i}^2 + \sum_{j=1}^{q} \beta_j \sigma_{t-j}^2, \qquad t \in \mathbb{Z},$$

where $(Z_t)$ is a sequence of i.i.d. random variables with $\mathbb{E}Z_0 = 0$ and $\mathbb{E}Z_0^2 = 1$, and $\alpha_i$ and $\beta_j$ are nonnegative parameters ensuring the nonnegativity of the squared *volatility* process $(\sigma_t^2)$. The GARCH model has attracted a lot of attention in the financial econometrics community because it constitutes an example of a *stationary* time series with a time varying *conditional* variance, hence, conditional heteroscedasticity. There is empirical evidence that GARCH captures some of the features of real-life financial return data rather well; see the survey article by Shephard [30] and the references therein.

Although the equations (1.1)–(1.2) look rather simple, they create a complicated probabilistic structure which is not easily understood. For example, the problem of finding conditions for the existence and uniqueness of a stationary solution to the equations (1.1)–(1.2) waited for a solution until Nelson [26] provided the answer for GARCH$(1, 1)$ and Bougerol and Picard [6] for the general GARCH$(p, q)$ case. Statistical inference leads to further nontrivial problems. Since the exact distribution of $Z_t$ remains unspecified, one determines likelihoods under the hypothesis of standard Gaussian innovations. Moreover, the volatility $\sigma_t$ is an *unobserved* quantity. A possible approximation of $\sigma_t^2$ is obtained by mimicking the recursion (1.2), for example, with an initialization $X_{-p+1} = \cdots = X_0 = 0$ and $\sigma_{-q+1}^2 = \cdots = \sigma_0^2 = 0$.

Recently, Berkes, Horváth and Kokoszka [2] showed under minimal assumptions that the resulting estimator, which is called the (Gaussian) quasi-maximum-likelihood estimator (QMLE), is consistent and asymptotically normal, thereby generalizing work by Lumsdaine [24], Lee and Hansen [22] and Boussama [7, 8]. Related to [2] is Francq and Zakoïan [15]. See also [1, 17, 25], where the QMLE in the GARCH$(p, q)$ model with a *heavy-tailed* noise distribution is studied.

Implicitly, the authors of [2] make use of the fact that GARCH$(p, q)$ can be embedded in a *stochastic recurrence equation* (SRE). We focus in this introduction on the particular GARCH$(1, 1)$ model in order to illustrate the method and to avoid heavy notation. The SRE for $\sigma_t^2$ is then given by

$$\sigma_t^2 = (\alpha_1 Z_{t-1}^2 + \beta_1)\sigma_{t-1}^2 + \alpha_0 = A_t \sigma_{t-1}^2 + B_t, \qquad t \in \mathbb{Z},$$

where $A_t = \alpha_1 Z_{t-1}^2 + \beta_1$ and $B_t = \alpha_0$; notice that $((A_t, B_t))$ constitutes an i.i.d. sequence. From the literature on SRE's (e.g., [9]), it is well known that the conditions $\mathbb{E}(\log|A_0|) < 0$ and $\mathbb{E}(\log^+|B_0|) < \infty$ guarantee the existence and uniqueness of a strictly stationary solution to the SRE

$$(1.3) \qquad Y_t = A_t Y_{t-1} + B_t, \qquad t \in \mathbb{Z},$$

provided $((A_t, B_t))$ is a stationary ergodic sequence. The condition $\mathbb{E}(\log|A_0|) < 0$ can be interpreted as "the $|A_t|$'s are smaller than one on average." We will refer to this property as the *contraction property* of the SRE. Altogether, applied to GARCH(1,1), this shows that $\mathbb{E}[\log(\alpha_1 Z_0^2 + \beta_1)] < 0$ is a sufficient condition for the existence of a stationary solution. Nelson [26] also gave an argument which shows that $\mathbb{E}[\log(\alpha_1 Z_0^2 + \beta_1)] < 0$ is necessary for the stationarity of GARCH(1,1). Note that the latter condition implies $\beta_1 < 1$ because $\log(\beta_1) \leq \mathbb{E}[\log(\alpha_1 Z_0^2 + \beta_1)] < 0$.

Now, suppose for the moment that the parameters $\alpha_0, \alpha_1, \beta_1$ of a stationary GARCH(1,1) are known. The fact that $\beta_1 < 1$ enables us to approximate the unobserved squared volatility through the recursion $\hat{\sigma}_t^2 = \alpha_0 + \alpha_1 X_{t-1}^2 + \beta_1 \hat{\sigma}_{t-1}^2$, $t \geq 1$, with initialization $\hat{\sigma}_0^2 \geq 0$. Indeed, the relation $|\hat{\sigma}_t^2 - \sigma_t^2| = (\beta_1)^t |\hat{\sigma}_0^2 - \sigma_0^2|$ shows that the error of approximation decays at an exponential rate when $t \to \infty$. Likewise, for an arbitrary parameter $\boldsymbol{\theta} = (u_1, u_2, u_3)^T$ with $0 \leq u_3 < 1$, one defines $\hat{h}_t(\boldsymbol{\theta}) = u_1 + u_2 X_{t-1}^2 + u_3 \hat{h}_{t-1}(\boldsymbol{\theta})$, $t \geq 1$, together with the initialization $\hat{h}_0(\boldsymbol{\theta}) \geq 0$. This means that the conditional log-likelihood of $(X_1, \ldots, X_n)^T$ given $(X_0, \sigma_0)^T$ under $Z_t$ i.i.d. $\sim \mathcal{N}(0,1)$ is approximately equal to

$$\hat{L}_n(\boldsymbol{\theta}) = -\frac{1}{2} \sum_{t=1}^{n} \left( \frac{X_t^2}{\hat{h}_t(\boldsymbol{\theta})} + \log \hat{h}_t(\boldsymbol{\theta}) \right),$$

modulo an unimportant constant. We define the QMLE as the maximizer $\hat{\boldsymbol{\theta}}_n$ of $\hat{L}_n$ with respect to $\boldsymbol{\theta}$. Because of $0 \leq u_3 < 1$, the sequence of random functions $(\hat{h}_t)_{t \geq 0}$ can be approximated by a stationary ergodic process $(h_t)$ such that the error $|\hat{h}_t - h_t|$ decays to zero uniformly on a certain compact set and with geometric rate. It turns out that $(h_t)$ is characterized as the unique stationary ergodic solution to the SRE

$$(1.4) \qquad h_t(\boldsymbol{\theta}) = u_1 + u_2 X_{t-1}^2 + u_3 h_{t-1}(\boldsymbol{\theta}), \qquad t \geq 1.$$

We will specify the relevant normed function space on which we study the SRE (1.4) in the subsequent sections. The authors of [2] use the uniform exponential decay of $|\hat{h}_t - h_t|$ in order to prove that the maximizer $\tilde{\boldsymbol{\theta}}_n$ of

$$L_n(\boldsymbol{\theta}) = -\frac{1}{2} \sum_{t=1}^{n} \left( \frac{X_t^2}{h_t(\boldsymbol{\theta})} + \log h_t(\boldsymbol{\theta}) \right)$$

is asymptotically equivalent to $\hat{\boldsymbol{\theta}}_n$. It is therefore enough to analyze $\tilde{\boldsymbol{\theta}}_n$. This is a slightly simpler problem because $((X_t^2/h_t + \log h_t))$ is a stationary ergodic sequence. The consistency of $\tilde{\boldsymbol{\theta}}_n$ then follows by an application of standard techniques, and the asymptotic normality of $\tilde{\boldsymbol{\theta}}_n$ is shown by a Taylor expansion of $L'_n$ together with an application of a martingale central limit theorem. Note that formal differentiation of (1.4) with respect to $\boldsymbol{\theta}$ leads to the SRE

$$(1.5) \qquad h'_t = (1, X^2_{t-1}, h_{t-1})^T + u_3 h'_{t-1}, \qquad t \geq 1.$$

Since $0 \leq u_3 < 1$, the linear SRE (1.5) is contractive and admits a unique stationary ergodic solution $(h'_t)_{t \geq 0}$. It can be shown that $h'_t$ indeed coincides with the first derivative of $h_t$ for all $t \geq 0$. By similar arguments, one establishes the existence, stationarity and ergodicity of the second derivatives of $h_t$, so that the Taylor expansion of $L'_n$ is validated.

The ideas sketched above can be extended and applied to conditionally heteroscedastic models of the more general form

$$(1.6) \qquad \begin{cases} X_t = \sigma_t Z_t, \\ \sigma^2_t = g_{\boldsymbol{\theta}}(X_{t-1}, \ldots, X_{t-p}, \sigma^2_{t-1}, \ldots, \sigma^2_{t-q}), \end{cases} \qquad t \in \mathbb{Z},$$

where $p, q \geq 0$ are integers, $\{g_{\boldsymbol{\theta}}\}$ denotes a parametric family of nonnegative functions on $\mathbb{R}^p \times [0, \infty)^q$ satisfying certain conditions and $(Z_t)$ is a sequence of i.i.d. random variables with $\mathbb{E}Z_0 = 0$ and $\mathbb{E}Z^2_0 = 1$. Essentially, this is accomplished under the assumption that "$g_{\boldsymbol{\theta}}$ is a contraction." This condition is a random coefficient Lipschitz condition, where the Lipschitz coefficient has a negative logarithmic moment, analogous to $\mathbb{E}(\log A_0) < 0$ in the GARCH(1, 1) case. The notion of contractivity will be made precise in different forms in the following sections. It is then possible to treat—in a unified manner—the QMLE for various other GARCH-type models, such as AGARCH (asymmetric GARCH), introduced by Ding, Granger and Engle [14], or EGARCH (exponential GARCH), introduced by Nelson [27].

The main goals of this paper are the following: (1) studying the general conditionally heteroscedastic model (1.6) with emphasis on the approximation $\hat{\sigma}_t$ of the stochastic volatility $\sigma_t$ given initial values and (2) quasi-maximum-likelihood inference. The novel approach in this paper is to *explicitly* formulate and solve the relevant problems by making use of stochastic recurrence equations (SRE's). This will lead to proofs which illuminate the underlying probabilistic structure.

It turns out that *invertibility* of the model, that is, the fact that

$$|\hat{\sigma}^2_t - \sigma^2_t| \xrightarrow{\mathbb{P}} 0, \qquad t \to \infty,$$

is an important condition for the applicability of the QMLE. The verification of the invertibility property is a rather straightforward task for the GARCH and AGARCH models, but establishing invertibility of EGARCH is an intricate and unpleasant problem. A final treatment of the QMLE in EGARCH is not possible at the time being, and one may regard this open problem as one of the limitations of this model. To the best of our knowledge, the theoretical properties of the QMLE in EGARCH have not been studied in the literature.

Work of Boussama [7, 8] and Comte and Lieberman [12] on multivariate GARCH indicates that the contraction technique cannot be employed in the analysis of the QMLE in multivariate conditionally heteroscedastic time series models.

The paper is organized as follows. In Section 2 we collect some of the technical tools to be used throughout this paper. Among them, Section 2.5 is crucial since it gives conditions for the existence and uniqueness of a strictly stationary solution to an SRE with values in abstract spaces. This section is based on work by Bougerol [5]. In Section 3 we investigate some of the basic properties of the general conditionally heteroscedastic model. We start with the problem of existence and uniqueness of a stationary solution to the underlying SRE (Section 3.1), discuss the invertibility of the model (Section 3.2) and finish with a proper definition of the above mentioned sequence $(h_t)$ in the general model (1.6). In Section 4 we establish the consistency of the QMLE in the general conditionally heteroscedastic model (1.6). In Section 5 we verify the conditions for consistency of the QMLE for two particular models, AGARCH$(p,q)$, including GARCH$(p,q)$, and EGARCH. In Section 6 we study the differentiability of $h_t$ as a preliminary but crucial step in the proof of the asymptotic normality of the QMLE in the general model (1.6) given in Section 7. We verify the conditions for asymptotic normality of AGARCH$(p,q)$ in Section 8. The verification of these conditions in the general EGARCH model is an open question. However, in some specific models asymptotic normality can be established as well. Since this leads to a series of technical calculations, we refer the reader to [31]. In Section 9 we briefly indicate that the QMLE is also applicable to certain nonstationary solutions of the model (1.6).

## 2. Preliminaries.

2.1. *Exponentially fast almost sure convergence.* Let $T = \mathbb{N} = \{0, 1, \dots\}$ or $T = \mathbb{Z}$. For sequences $(v_t)_{t \in \mathbb{Z}}$, we write $(v_t)$ in abridged form. The following definition will facilitate the formulation of many of our results. A sequence $(v_t)_{t \in T}$ of random elements with values in a normed vector space $(B, \|\cdot\|)$ is said to converge to zero *exponentially fast* a.s. when $t \to \infty$, henceforth, denoted by $v_t \overset{\text{e.a.s.}}{\longrightarrow} 0$, if there exists $\gamma > 1$ with $\gamma^t \|v_t\| \overset{\text{a.s.}}{\longrightarrow} 0$.

LEMMA 2.1. *Let $(\xi_t)_{t \in T}$ be a sequence of real random variables with $\xi_t \overset{\text{e.a.s.}}{\longrightarrow} 0$ and $(v_t)_{t \in T}$ a sequence of identically distributed random elements* *with values in a separable Banach space* $(B, \|\cdot\|)$. *If* $\mathbb{E}(\log^+ \|v_0\|) < \infty$, *then* $\sum_{t=0}^\infty \xi_t v_t$ *converges a.s., and one has* $\xi_n \sum_{t=0}^n v_t \xrightarrow{\text{e.a.s.}} 0$ *and* $\xi_n v_n \xrightarrow{\text{e.a.s.}} 0$ *as* $n \to \infty$.

PROOF.   Since $\xi_t \xrightarrow{\text{e.a.s.}} 0$, there is a $\gamma > 1$ with $\gamma^t |\xi_t| \xrightarrow{\text{a.s.}} 0$. Hence, there exists a random variable $C \geq 0$ with $|\xi_t| \leq C\gamma^{-t}$ a.s. for all $t \in \mathbb{N}$. It is a well-known fact that $\gamma^{-1} < 1$ and $\mathbb{E}(\log^+ \|v_0\|) < \infty$ imply $\sum_{t=0}^\infty \gamma^{-t}\|v_t\| < \infty$ a.s.; see, for example, [2]. Therefore,

$$\sum_{t=0}^\infty \|\xi_t v_t\| = \sum_{t=0}^\infty |\xi_t| \|v_t\| \leq C \sum_{t=0}^\infty \gamma^{-t}\|v_t\| < \infty \qquad \text{a.s.}$$

Since $B$ is complete, this implies that $\sum_{t=0}^\infty \xi_t v_t$ converges in $B$ a.s. As to the second part of the lemma, choose $\tilde\gamma > 1$ small enough such that $\eta = \tilde\gamma/\gamma < 1$. Since $\eta^{1/2} < 1$ and $\mathbb{E}(\log^+ \|v_0\|) < \infty$ imply $\sum_{t=0}^\infty \eta^{t/2}\|v_t\| < \infty$ a.s., the estimate

$$\tilde\gamma^n \left\| \xi_n \sum_{t=0}^n v_t \right\| = \eta^n \gamma^n |\xi_n| \left\| \sum_{t=0}^n v_t \right\| \leq C\eta^{n/2} \sum_{t=0}^n \eta^{n/2}\|v_t\| \leq C\eta^{n/2} \sum_{t=0}^\infty \eta^{t/2}\|v_t\|$$

establishes $\xi_n \sum_{t=0}^n v_t \xrightarrow{\text{e.a.s.}} 0$ as $n \to \infty$ and $\xi_n v_n \xrightarrow{\text{e.a.s.}} 0$. This completes the proof.   □

We will often use the following elementary properties without any explicit reference.

LEMMA 2.2.   *Let* $v_1$ *and* $v_2$ *be two random elements taking values in a separable Banach space. Then the existence of a* $q > 0$ *such that* $\mathbb{E}\|v_1\|^q < \infty$ *implies* $\mathbb{E}(\log^+ \|v_1\|) < \infty$. *Moreover,*

$$\mathbb{E}(\log^+ \|v_1 + v_2\|) \leq 2\log 2 + \mathbb{E}(\log^+ \|v_1\|) + \mathbb{E}(\log^+ \|v_2\|)$$

*and*

$$\mathbb{E}[\log^+(\|v_1\|\|v_2\|)] \leq \mathbb{E}(\log^+ \|v_1\|) + \mathbb{E}(\log^+ \|v_2\|).$$

We also need the following auxiliary result.

LEMMA 2.3.   *Let* $(Y_t)_{t \in T}$ *and* $(\tilde Y_t)_{t \in T}$ *be sequences with* $|Y_t - \tilde Y_t| \xrightarrow{\text{e.a.s.}} 0$ *as* $t \to \infty$. *Then we have the following:*

   (i)  $|Y_t|^{1/2} - |\tilde Y_t|^{1/2} \xrightarrow{\text{e.a.s.}} 0$ *as* $t \to \infty$.
   (ii)  *If, in addition,* $(Y_t)_{t \in T}$ *is stationary with* $\mathbb{E}|Y_0| < \infty$, *then also* $|\exp(Y_t) - \exp(\tilde Y_t)| \xrightarrow{\text{e.a.s.}} 0$.

PROOF.    (i) By assumption on $|Y_t - \tilde{Y}_t|$, there is a $\gamma > 1$ with $\gamma^t |Y_t - \tilde{Y}_t| \xrightarrow{\text{a.s.}} 0$. Consequently, there exists $t_0$ (depending on the realization $\omega$) such that, for every $t \geq t_0$, either $\{\min(Y_t, \tilde{Y}_t) \leq \gamma^{-t}, \max(Y_t, \tilde{Y}_t) \leq \gamma^{-t/2}\}$ or $\{\min(Y_t, \tilde{Y}_t) > \gamma^{-t}\}$. Using this observation and applying the mean value theorem to the difference of square roots, we conclude

$$|(Y_t)^{1/2} - (\tilde{Y}_t)^{1/2}| \leq \max\{\gamma^{-t/2}, 2^{-1}\gamma^{t/2}|Y_t - \tilde{Y}_t|\}, \qquad t \geq t_0.$$

This shows that $\tilde{\gamma}^t |(Y_t)^{1/2} - (\tilde{Y}_t)^{1/2}| \xrightarrow{\text{a.s.}} 0$ if $1 < \tilde{\gamma} < \gamma^{1/2}$ and, thus, $|(Y_t)^{1/2} - (\tilde{Y}_t)^{1/2}| \xrightarrow{\text{e.a.s.}} 0$.

(ii) As soon as $|Y_t - \tilde{Y}_t| \leq 1$, by the mean value theorem,

$$|\exp(Y_t) - \exp(\tilde{Y}_t)| = \exp(Y_t)|1 - \exp(\tilde{Y}_t - Y_t)| \leq \exp(1 + Y_t)|Y_t - \tilde{Y}_t|.$$

Now the claim follows from $|Y_t - \tilde{Y}_t| \xrightarrow{\text{e.a.s.}} 0$ together with $\mathbb{E}[\log^+\{\exp(1 + Y_0)\}] \leq 1 + \mathbb{E}|Y_0| < \infty$ and an application of Lemma 2.1. This concludes the proof.   □

The following elementary result about random products is used in abundance.

LEMMA 2.4.    Let $(Y_t)_{t \in T}$ be a stationary ergodic sequence of nonnegative random variables with $\mathbb{E}(\log Y_0) < 0$. Then

$$\prod_{t=0}^{n} Y_t \xrightarrow{\text{e.a.s.}} 0, \qquad n \to \infty.$$

If, in addition, $(Y_t)_{t \in T}$ is an i.i.d. sequence with $\mathbb{E}Y_0^q < \infty$ for some $q > 0$, then there exist $0 < \tilde{q} \leq q$ and $0 < \eta < 1$ such that

$$\mathbb{E}\left(\prod_{t=0}^{n-1} Y_t\right)^{\tilde{q}} = \eta^n, \qquad n \to \infty.$$

PROOF.    The first assertion follows from a straightforward application of the ergodic theorem to the logarithm of $\prod_{t=0}^{n} Y_t$ together with $\mathbb{E}[\log Y_0] < 0$. Concerning the second claim, observe that the map $s \mapsto \mathbb{E}|Y_0|^s$ on $[0, q]$ has first derivative equal to $\mathbb{E}(\log Y_0) < 0$ at $s = 0$. This shows the existence of $\tilde{q} \in (0, q]$ with $\eta := \mathbb{E}Y_0^{\tilde{q}} < 1$, and since $(Y_t)_{t \in T}$ is i.i.d., $\mathbb{E}(\prod_{t=0}^{n-1} Y_t)^{\tilde{q}} = (\mathbb{E}Y_0^{\tilde{q}})^n = \eta^n$, which concludes the proof of the lemma.   □

**2.2. Stationarity and ergodicity.**    The proofs of the consistency and asymptotic normality of the QMLE rest on the ergodicity of the underlying time series. The stationarity and ergodicity of sequences $(v_t)_{t \in T}$ of random elements with values in a general measurable space are up to some evident generalizations defined as in the case of sequences of real random variables; see, for example, the monograph of Krengel [20]. A standard example of a stationary ergodic sequence of random elements is an i.i.d. sequence. For the following well-known result, see Proposition 4.3 in [20].

PROPOSITION 2.5. *Let* $(E, \mathcal{E})$ *and* $(\tilde{E}, \tilde{\mathcal{E}})$ *be measurable spaces. If* $(v_t)$ *is a stationary ergodic sequence of* $E$-*valued random elements and* $h\colon (E^{\mathbb{N}}, \mathcal{E}^{\mathbb{N}}) \to (\tilde{E}, \tilde{\mathcal{E}})$ *is measurable, then the sequence* $(\tilde{v}_t)$ *defined by*

$$(2.1) \qquad\qquad \tilde{v}_t = h(v_t, v_{t-1}, \dots), \qquad t \in \mathbb{Z},$$

*is stationary ergodic.*

Since we did not find a proof for the following result in the literature, we refer to [31] for details.

PROPOSITION 2.6. *Let* $(E, \mathcal{E})$ *be a measurable space and* $\tilde{E}$ *a complete and separable metric space, which we endow with its Borel* $\sigma$-*field* $\tilde{\mathcal{E}} = \mathcal{B}(\tilde{E})$. *Assume that* $(v_t)$ *is a stationary ergodic sequence of* $E$-*valued random elements. Let* $f_m\colon (E^{\mathbb{N}}, \mathcal{E}^{\mathbb{N}}) \to (\tilde{E}, \tilde{\mathcal{E}})$, $m \in \mathbb{N}$, *be measurable maps such that for some* $t_0 \in \mathbb{Z}$, *the sequence* $f_m(v_{t_0-1}, v_{t_0-2}, \dots)$ *is convergent in* $\tilde{E}$ *a.s. Then there is a measurable map* $f\colon (E^{\mathbb{N}}, \mathcal{E}^{\mathbb{N}}) \to (\tilde{E}, \tilde{\mathcal{E}})$ *such that, for all* $t \in \mathbb{Z}$,

$$\tilde{v}_t = \lim_{m\to\infty} f_m(v_{t-1}, v_{t-2}, \dots) = f(v_{t-1}, v_{t-2}, \dots) \qquad a.s.,$$

*and the sequence* $(\tilde{v}_t)$ *is stationary ergodic.*

2.3. *Uniform convergence via the ergodic theorem.* For establishing consistency and asymptotic normality of $M$-estimators, one often needs to show the uniform convergence of sequences of random functions. Let us assume that $K \subset \mathbb{R}^d$ is a compact set. Write $\mathbb{C}(K, \mathbb{R}^{d'})$ for the space of continuous $\mathbb{R}^{d'}$-valued functions equipped with the sup-norm $\|v\|_K = \sup_{s \in K} |v(s)|$, and for short, $\mathbb{C}(K, \mathbb{R}) = \mathbb{C}(K)$. Recall that $\mathbb{C}(K, \mathbb{R}^{d'})$ is a separable Banach space. For a given stationary sequence of random elements $(v_t)_{t \in T}$ with values in $\mathbb{C}(K, \mathbb{R}^{d'})$, we say that the *uniform strong law of large numbers* holds, if in $\mathbb{C}(K, \mathbb{R}^{d'})$,

$$\frac{1}{n}\sum_{t=1}^{n} v_t \xrightarrow{\text{a.s.}} v, \qquad n \to \infty,$$

where the function $v$ is defined pointwise by $v(s) = \mathbb{E}[v_0(s)]$, $s \in K$. Note that the a.s. uniform convergence understands that $v \in \mathbb{C}(K, \mathbb{R}^{d'})$. In the literature the uniform SLLN is often established in two steps. First, by an application of the ergodic theorem for *real* sequences, $n^{-1} \sum_{t=1}^{n} v_t(s) \xrightarrow{\text{a.s.}} v(s)$, for every *fixed* $s \in K$. Second, the uniform convergence is proved by showing the a.s. equicontinuity of $\{n^{-1} \sum_{t=1}^{n} v_t\}$. In order to establish the a.s. equicontinuity of the *second* derivatives of the log-likelihood, Lumsdaine [24], Lee and Hansen [22] and Berkes, Horváth and Kokoszka [2] stochastically bound the *third* derivatives. Such computations can be avoided when the ergodic theorem for random elements with values in a separable Banach space is applied. We state the relevant result as a theorem and refer to [29] for its proof.

THEOREM 2.7. *Let* $(v_t)$ *be a stationary ergodic sequence of random elements with values in* $\mathbb{C}(K, \mathbb{R}^{d'})$. *Then the uniform SLLN is implied by* $\mathbb{E}\|v_0\|_K < \infty$.

2.4. *Notation and matrix norms.* For a sequence $(\phi_t)$ of transformations on a certain space $E$, we denote by $(\phi_t^{(r)})$ the sequence of the $n$-fold iterations of past and present transformations defined by

$$\phi_t^{(r)} = \begin{cases} \mathrm{Id}_E, & r = 0, \\ \phi_t \circ \phi_{t-1} \circ \cdots \circ \phi_{t-r+1}, & r \geq 1, \end{cases}$$

where $\mathrm{Id}_E$ denotes the identity map in $E$.

In this paper we use two different matrix norms, as the case may be; since all matrix norms are equivalent, the particular choice of a norm is (mostly) irrelevant from a mathematical point of view. Recall that the *Frobenius* norm of a matrix $\mathbf{A} = (a_{ij}) \in \mathbb{R}^{d' \times d'}$ is defined by

$$\|\mathbf{A}\| = \left( \sum_{i,j} a_{ij}^2 \right)^{1/2}.$$

Occasionally we use the *operator norm* with respect to the Euclidean norm instead, that is,

$$\|\mathbf{A}\|_{\mathrm{op}} = \sup_{\mathbf{x} \neq \mathbf{0}} \frac{|\mathbf{A}\mathbf{x}|}{|\mathbf{x}|}.$$

Finally we define the norm of a continuous matrix-valued function $\mathbf{A}$ on a compact set $K \subset \mathbb{R}^d$, that is, $\mathbf{A} \in \mathbb{C}(K, \mathbb{R}^{d' \times d'})$, by

$$\|\mathbf{A}\|_K = \sup_{\mathbf{s} \in K} \|\mathbf{A}(\mathbf{s})\|.$$

Of course, Theorem 2.7 carries over to sequences of matrices: a stationary ergodic sequence of random elements $(\mathbf{A}_t)$ with values in $\mathbb{C}(K, \mathbb{R}^{d' \times d'})$ and with $\mathbb{E}\|\mathbf{A}_0\|_K < \infty$ satisfies

$$\frac{1}{n} \sum_{t=1}^{n} \mathbf{A}_t \xrightarrow{\text{a.s.}} \mathbf{M}_0, \qquad n \to \infty,$$

where $\mathbf{M}_0(\mathbf{s}) = \mathbb{E}[\mathbf{A}_0(\mathbf{s})]$, $\mathbf{s} \in K$; note that the expected value of a random matrix is defined element-wise.

2.5. *Stochastic recurrence equations.* Although we merely work in Banach spaces, we present the notion of an SRE in the more general case of a complete separable metric space (i.e., Polish space), as it was, for example, formulated by Bougerol [5]. Hence, let $(E, d)$ be a Polish space equipped with its Borel $\sigma$-field $\mathcal{B}(E)$. Recall that a map $\phi : E \to E$ is called *Lipschitz* if

$$(2.2) \qquad \Lambda(\phi) := \sup_{x, y \in E, x \neq y} \left( \frac{d(\phi(x), \phi(y))}{d(x, y)} \right)$$

is finite and is called a *contraction* if $\Lambda(\phi) < 1$. Also note that $\Lambda$ is submultiplicative, that is, if $\phi$ and $\psi$ are Lipschitz maps $E \to E$, then

$$(2.3) \qquad \Lambda(\phi \circ \psi) \leq \Lambda(\phi) \Lambda(\psi).$$

We consider a process $(\phi_t)$ of random Lipschitz maps $E \to E$ with $\phi_t(x)$ being $\mathcal{B}(E)$-measurable for every fixed $x \in E$ and $t \in \mathbb{Z}$. In what follows, $T = \mathbb{N}$ or $T = \mathbb{Z}$. If for a stochastic process $(X_t)_{t \in T}$ with values in $E$,

$$(2.4) \qquad X_{t+1} = \phi_t(X_t), \qquad t \in T,$$

we say that $(X_t)_{t \in T}$ obeys the SRE associated with $(\phi_t)$. Alternatively, $(X_t)_{t \in T}$ is referred to as a *solution* to the SRE (2.4). The following theorem due to Bougerol [5] is a stochastic version of Banach's fixed point theorem and makes statements about the solutions of (2.4) under the assumption that $(\phi_t)$ is stationary ergodic and that $\phi_0$ or a certain $r$-fold iterate $\phi_0^{(r)}$ is "contractive on average."

THEOREM 2.8 (Theorem 3.1 of [5]). *Let $(\phi_t)$ be a stationary ergodic sequence of Lipschitz maps from $E$ into $E$. Suppose the following conditions hold:*

S.1 *There is a $y \in E$ such that $\mathbb{E}[\log^+ d(\phi_0(y), y)] < \infty$.*
S.2 *$\mathbb{E}[\log^+ \Lambda(\phi_0)] < \infty$ and for some integer $r \geq 1$,*

$$\mathbb{E}[\log \Lambda(\phi_0^{(r)})] = \mathbb{E}[\log \Lambda(\phi_0 \circ \cdots \circ \phi_{-r+1})] < 0.$$

*Then the SRE (2.4) admits a stationary solution $(Y_t)_{t \in T}$ which is ergodic. The following almost sure representation is valid:*

$$(2.5) \qquad Y_t = \lim_{m \to \infty} \phi_{t-1} \circ \cdots \circ \phi_{t-m}(y), \qquad t \in T,$$

*where the limit is irrespective of $y \in E$. The random elements $Y_t$ are measurable with respect to the $\sigma$-field generated by $\{\phi_{t-k} \mid k \geq 1\}$. If $(\tilde{Y}_t)_{t \in T}$ is any other solution to (2.4), then*

$$(2.6) \qquad d(\tilde{Y}_t, Y_t) \xrightarrow{\text{e.a.s.}} 0, \qquad t \to \infty.$$

*Moreover, in the case $T = \mathbb{Z}$ the stationary solution to the SRE (2.4) is unique.*

REMARK 2.9. 1. Note that there are in general many solutions to (2.4). Indeed, if $T = \mathbb{N}$, take any $z \in E$ in order to produce a solution $(\phi_{t-1}^{(t)}(z))_{t \in \mathbb{N}}$. The elements $\phi_{t-1}^{(t)}(z) = \phi_{t-1} \circ \cdots \circ \phi_0(z)$ are also called the *forward iterates* associated with the SRE (2.4), whereas the elements $\phi_{t-1} \circ \cdots \circ \phi_{t-m}(z)$ with $t$ fixed and $m \geq 0$ are called *backward iterates*. Thus, property (2.5) can be read as $Y_t$ being the limit of its backward iterates, and relation (2.6) means that the forward iterates approach the trajectory of the unique stationary solution $(Y_t)_{t \in T}$ with an error decaying exponentially fast as $t \to \infty$.

2. Although it is an easy matter to demonstrate the uniqueness of the stationary solution $(Y_t)$, this property is not mentioned in Theorem 3.1 of [5]. We give an argument for the sake of completeness. Indeed, if $(\tilde{Y}_t)$ is yet another stationary solution of the SRE (2.4), then by the definition (2.2) and the triangle inequality,

$$
\begin{aligned}
(2.7) \quad d(\tilde{Y}_t, Y_t) &= d(\phi_{t-1}^{(m)}(\tilde{Y}_{t-m}), \phi_{t-1}^{(m)}(Y_{t-m})) \leq \Lambda(\phi_{t-1}^{(m)})\, d(\tilde{Y}_{t-m}, Y_{t-m}) \\
&\leq \Lambda(\phi_{t-1}^{(m)})(d(\tilde{Y}_{t-m}, y) + d(y, Y_{t-m})).
\end{aligned}
$$

As shown in [5], as a consequence of S.2, one has

$$
(2.8) \qquad\qquad \Lambda(\phi_{t-1}^{(m)}) \xrightarrow{\text{e.a.s.}} 0, \qquad m \to \infty.
$$

Since $(d(\tilde{Y}_{t-m}, y))_{m \in \mathbb{N}}$ and $(d(y, Y_{t-m}))_{m \in \mathbb{N}}$ are stationary, a Slutsky argument shows that the right-hand side of (2.7) converges to zero in probability. For this reason, $\mathbb{P}(d(\tilde{Y}_t, Y_t) = 0) = 1$, and $(\tilde{Y}_t)_{t \in T}$ and $(Y_t)_{t \in T}$ are indistinguishable.

The following result relates the solutions of an SRE associated with a stationary ergodic sequence $(\phi_t)$ of Lipschitz maps with the solutions of a certain perturbed SRE.

THEOREM 2.10. *Let $(B, \|\cdot\|)$ be a separable Banach space and $(\phi_t)$ be a stationary ergodic sequence of Lipschitz maps from $B$ into $B$. Impose the following:*

S.1 $\mathbb{E}(\log^+ \|\phi_0(0)\|) < \infty$.
S.2 $\mathbb{E}[\log^+ \Lambda(\phi_0)] < \infty$ *and for some integer $r \geq 1$,*

$$
\mathbb{E}[\log \Lambda(\phi_0^{(r)})] = \mathbb{E}[\log \Lambda(\phi_0 \circ \cdots \circ \phi_{-r+1})] < 0.
$$

*Assume that $\mathbb{E}(\log^+ \|Y_0\|) < \infty$ for the stationary ergodic solution $(Y_t)_{t \in \mathbb{N}}$ of the SRE associated with $(\phi_t)$, which is given by (2.5). Let $(\hat{\phi}_t)_{t \in \mathbb{N}}$ be a sequence of Lipschitz maps such that*

S.3 $\|\hat{\phi}_t(0) - \phi_t(0)\| \xrightarrow{\text{e.a.s.}} 0$ *and $\Lambda(\hat{\phi}_t - \phi_t) \xrightarrow{\text{e.a.s.}} 0$ as $t \to \infty$.*

*Then for every solution* $(\hat{Y}_t)_{t\in\mathbb{N}}$ *of the perturbed SRE,*

$$X_{t+1} = \hat{\phi}_t(X_t), \qquad t \in \mathbb{N},$$

*one has that*

$$(2.9) \qquad \qquad \|\hat{Y}_t - Y_t\| \overset{\text{e.a.s.}}{\longrightarrow} 0, \qquad t \to \infty.$$

PROOF. Note that the sequence $(\phi_t)$ satisfies the conditions of Theorem 2.8 with $d$ the metric induced by $\|\cdot\|$. It is sufficient to demonstrate $\|\hat{Y}_{s+kr} - Y_{s+kr}\| \overset{\text{e.a.s.}}{\longrightarrow} 0$ as $k \to \infty$ for each $s \in [0, r)$; indeed, the latter limit relation implies $\|\hat{Y}_t - Y_t\| \le \sum_{s=0}^{r-1} \|\hat{Y}_{s+r\lfloor t/r\rfloor} - Y_{s+r\lfloor t/r\rfloor}\| \overset{\text{e.a.s.}}{\longrightarrow} 0$ as $t \to \infty$. To begin with, we establish the auxiliary results

$$(2.10) \qquad \begin{aligned} d_t &:= \|\hat{\phi}_t^{(r)}(0) - \phi_t^{(r)}(0)\| \overset{\text{e.a.s.}}{\longrightarrow} 0 \quad \text{and} \\ e_t &:= \Lambda(\hat{\phi}_t^{(r)} - \phi_t^{(r)}) \overset{\text{e.a.s.}}{\longrightarrow} 0, \qquad t \to \infty, \end{aligned}$$

that is, the conditions of S.3 are also satisfied for the $r$-fold convolutions $\hat{\phi}_t^{(r)}$ and $\phi_t^{(r)}$. Note that the limit relation (2.10) is true if $r = 1$ by virtue of S.3. The proof of (2.10) proceeds by induction on $r$. Indeed, by the triangle inequality, for any $m \ge 1$,

$$\begin{aligned} &\|\hat{\phi}_t^{(m+1)}(0) - \phi_t^{(m+1)}(0)\| \\ &\quad = \|\hat{\phi}_t \circ \hat{\phi}_{t-1}^{(m)}(0) - \phi_t \circ \phi_{t-1}^{(m)}(0)\| \\ &\quad \le \|\hat{\phi}_t \circ \hat{\phi}_{t-1}^{(m)}(0) - \hat{\phi}_t \circ \phi_{t-1}^{(m)}(0)\| \\ &\qquad + \|(\hat{\phi}_t - \phi_t) \circ \phi_{t-1}^{(m)}(0) - (\hat{\phi}_t - \phi_t)(0)\| + \|(\hat{\phi}_t - \phi_t)(0)\| \\ &\quad \le \Lambda(\hat{\phi}_t)\|\hat{\phi}_{t-1}^{(m)}(0) - \phi_{t-1}^{(m)}(0)\| + \Lambda(\hat{\phi}_t - \phi_t)\|\phi_{t-1}^{(m)}(0)\| + \|\hat{\phi}_t(0) - \phi_t(0)\|. \end{aligned}$$

The e.a.s. convergence to zero of the left-hand side is a consequence of $\Lambda(\hat{\phi}_t) \le \Lambda(\hat{\phi}_t - \phi_t) + \Lambda(\phi_t)$, $\mathbb{E}[\log^+ \Lambda(\phi_0)] < \infty$, the limit relation $\|\hat{\phi}_{t-1}^{(\ell)}(0) - \phi_{t-1}^{(\ell)}(0)\| \overset{\text{e.a.s.}}{\longrightarrow} 0$ for every $\ell \in [0, m]$, and $\mathbb{E}(\log^+ \|\phi_{t-1}^{(m)}(0)\|) < \infty$ together with repeated application of Lemma 2.1. Exploiting the submultiplicativity of $\Lambda$, one finds in a similar way that

$$\begin{aligned} &\Lambda(\hat{\phi}_t^{(m+1)} - \phi_t^{(m+1)}) \\ &\quad \le \Lambda(\hat{\phi}_t)\Lambda(\hat{\phi}_{t-1}^{(m)} - \phi_{t-1}^{(m)}) + \Lambda(\hat{\phi}_t - \phi_t)\Lambda(\phi_{t-1}^{(m)}) \overset{\text{e.a.s.}}{\longrightarrow} 0, \qquad t \to \infty, \end{aligned}$$

and thus, (2.10) is established.

If we set $\hat{c}_t = \Lambda(\hat{\phi}_t^{(r)})$ and take into account that any map $\phi$ satisfies $\|\phi(x)\| \le \|\phi(x) - \phi(0)\| + \|\phi(0)\| \le \Lambda(\phi)\|x\| + \|\phi(0)\|$, $x \in B$, we obtain

$$\|\hat{Y}_{s+kr} - Y_{s+kr}\|$$

$$= \|\hat{\phi}_{s+kr-1}^{(r)}(\hat{Y}_{s+(k-1)r}) - \phi_{s+kr-1}^{(r)}(Y_{s+(k-1)r})\|$$

$$\leq \|\hat{\phi}_{s+kr-1}^{(r)}(\hat{Y}_{s+(k-1)r}) - \hat{\phi}_{s+kr-1}^{(r)}(Y_{s+(k-1)r})\|$$

$$\qquad + \|\hat{\phi}_{s+kr-1}^{(r)}(Y_{s+(k-1)r}) - \phi_{s+kr-1}^{(r)}(Y_{s+(k-1)r})\|$$

$$\leq \Lambda(\hat{\phi}_{s+kr-1}^{(r)})\|\hat{Y}_{s+(k-1)r} - Y_{s+(k-1)r}\| + \Lambda(\hat{\phi}_{s+kr-1}^{(r)} - \phi_{s+kr-1}^{(r)})\|Y_{s+(k-1)r}\|$$

$$\qquad + \|\hat{\phi}_{s+kr-1}^{(r)}(0) - \phi_{s+kr-1}^{(r)}(0)\|$$

$$= \hat{c}_{s+kr-1}\|\hat{Y}_{s+(k-1)r} - Y_{s+(k-1)r}\| + (e_{s+kr-1}\|Y_{s+(k-1)r}\| + d_{s+kr-1}).$$

Iterating the latter inequality until $k = 1$, we receive the final estimate

$$
\begin{aligned}
(2.11) \qquad & \|\hat{Y}_{s+kr} - Y_{s+kr}\| \\
& \leq \left(\prod_{\ell=1}^{k} \hat{c}_{s+\ell r-1}\right)\|\hat{Y}_s - Y_s\| \\
& \quad + \sum_{\ell=1}^{k}\left(\prod_{i=\ell+1}^{k} \hat{c}_{s+ir-1}\right)(e_{s+\ell r-1}\|Y_{s+(\ell-1)r}\| + d_{s+\ell r-1}).
\end{aligned}
$$

We show that the right-hand side of the latter inequality tends to zero e.a.s. Set $c_t = \Lambda(\phi_t^{(r)})$. An application of the monotone convergence theorem to $\mathbb{E}[\log^-(c_0 + \varepsilon)]$ and an application of the dominated convergence theorem to $\mathbb{E}[\log^+(c_0 + \varepsilon)]$ show $\mathbb{E}[\log(c_0 + \varepsilon)] \to \mathbb{E}[\log c_0]$ as $\varepsilon \downarrow 0$. Thus, there is $\varepsilon_0 > 0$ such that $\mathbb{E}[\log(c_0 + \varepsilon_0)] < 0$, and from Lemma 2.4, it follows that $\prod_{\ell=1}^{k}(c_{s+\ell r-1} + \varepsilon_0) \xrightarrow{\text{e.a.s.}} 0$ as $k \to \infty$. This limit relation together with the facts that $\hat{c}_t = \Lambda((\hat{\phi}_t^{(r)} - \phi_t^{(r)}) + \phi_t^{(r)}) \leq e_t + c_t$ and a.s. $e_t \leq \varepsilon_0$ for all but finitely many $t$'s leads to

$$(2.12) \qquad \prod_{\ell=1}^{k} \hat{c}_{s+\ell r-1} \leq \prod_{\ell=1}^{k}(e_{s+\ell r-1} + c_{s+\ell r-1}) \xrightarrow{\text{e.a.s.}} 0, \qquad k \to \infty,$$

and shows that $(\prod_{\ell=1}^{k} \hat{c}_{s+\ell r-1})\|\hat{Y}_s - Y_s\| \xrightarrow{\text{e.a.s.}} 0$. As regards the second term in (2.11), we use Lemma 2.1 together with (2.10) and $\mathbb{E}[\log^+\|Y_0\|] < \infty$ to prove that $p_{s+\ell r-1} := e_{s+\ell r-1}\|Y_{s+(\ell-1)r}\| + d_{s+\ell r-1} \xrightarrow{\text{e.a.s.}} 0$ when $\ell \to \infty$. This limit result and (2.12) imply the existence of a random variable $a_0$ and a number $0 < \gamma < 1$ so that $p_{s+\ell r-1} \leq a_0 \gamma^\ell$ for all $\ell \geq 1$. From this bound, we obtain

$$\sum_{\ell=1}^{k}\left(\prod_{i=\ell+1}^{k} \hat{c}_{s+ir-1}\right)p_{s+\ell r-1} \leq a_0 \sum_{\ell=1}^{k}\left(\prod_{i=\ell+1}^{k} \hat{c}_{s+ir-1}\right)\gamma^\ell.$$

Since the map $\delta \mapsto \mathbb{E}[\log(c_0 + \delta)]$ is continuous and increasing and $\mathbb{E}(\log c_0) < 0$, there exists a $\delta_0 > 0$ such that $0 > \mathbb{E}[\log(c_0 + \delta_0)] > \log \gamma$. Set $\tilde{\hat{c}}_t = \hat{c}_t + \delta_0$ and $\tilde{c}_t = c_t + \delta_0$, and note that we are done if we can show $\sum_{\ell=1}^{k} (\prod_{i=\ell+1}^{k} \hat{\tilde{c}}_{s+ir-1}) \gamma^{\ell} \xrightarrow{\text{e.a.s.}} 0$. Since $(\tilde{c}_t)$ is stationary ergodic and $\mathbb{E}(\log \tilde{c}_0) > \log \gamma$, we have that $\gamma^{\ell} (\prod_{i=1}^{\ell} \tilde{c}_{s+ir-1})^{-1} \xrightarrow{\text{e.a.s.}} 0$. Using the mean value theorem together with $\hat{\tilde{c}}_t, \tilde{c}_t \geq \delta_0$, we conclude

$$|(\hat{\tilde{c}}_{s+ir-1})^{-1} - (\tilde{c}_{s+ir-1})^{-1}| \leq (1/\delta_0^2)|\hat{\tilde{c}}_{s+ir-1} - \tilde{c}_{s+ir-1}| \xrightarrow{\text{e.a.s.}} 0, \qquad i \to \infty,$$

and, therefore, the identical arguments as used for (2.12) yield that also $\gamma^{\ell} (\prod_{i=1}^{\ell} \hat{\tilde{c}}_{s+ir-1})^{-1} \xrightarrow{\text{e.a.s.}} 0$. This implies the existence of a random variable $b_0$ with $\gamma^{\ell} (\prod_{i=1}^{\ell} \hat{\tilde{c}}_{s+ir-1})^{-1} \leq b_0$ for all $\ell \geq 1$, so that

$$\sum_{\ell=1}^{k} \left( \prod_{i=\ell+1}^{k} \hat{\tilde{c}}_{s+ir-1} \right) \gamma^{\ell} \leq b_0 k \prod_{i=1}^{k} \hat{\tilde{c}}_{s+ir-1} \xrightarrow{\text{e.a.s.}} 0, \qquad k \to \infty.$$

For the last step, we have used $\prod_{i=1}^{k} \hat{\tilde{c}}_{s+ir-1} \xrightarrow{\text{e.a.s.}} 0$, which is a consequence of $\mathbb{E}(\log \tilde{c}_0) < 0$ and $|\hat{\tilde{c}}_t - \tilde{c}_t| \xrightarrow{\text{e.a.s.}} 0$ together with the same arguments as applied in the derivation of (2.12). This completes the proof.  $\square$

## 3. The model: stationarity, ergodicity and invertibility.

We consider the general conditionally heteroscedastic model

$$(3.1) \quad \begin{cases} X_t = \sigma_t Z_t, \\ \sigma_t^2 = g_{\boldsymbol{\theta}}(X_{t-1}, X_{t-2}, \ldots, X_{t-p}; \sigma_{t-1}^2, \sigma_{t-2}^2, \ldots, \sigma_{t-q}^2), \end{cases} \qquad t \in \mathbb{Z},$$

where $\{g_{\boldsymbol{\theta}} : \boldsymbol{\theta} \in \boldsymbol{\theta}\}$ denotes a parametric family of nonnegative functions on $\mathbb{R}^p \times [0, \infty)^q$ and $(Z_t)$ is a sequence of i.i.d. random variables with $\mathbb{E}Z_0 = 0$ and $\mathbb{E}Z_0^2 = 1$. We also assume that $\sigma_t$ is nonnegative and $\mathcal{F}_{t-1}$-measurable, where

$$\mathcal{F}_t = \sigma(Z_k; k \leq t), \qquad t \in \mathbb{Z},$$

denotes the $\sigma$-field generated by the random variables $\{Z_k; k \leq t\}$. In the context of the latter model we introduce the notation

$$\mathbf{X}_t = (X_t, \ldots, X_{t-p+1})^T \quad \text{and} \quad \boldsymbol{\sigma}_t^2 = (\sigma_t^2, \ldots, \sigma_{t-q+1}^2)^T,$$

and write $\sigma_t^2 = g_{\boldsymbol{\theta}}(\mathbf{X}_{t-1}, \boldsymbol{\sigma}_{t-1}^2)$ for short.

In Section 3.1 we investigate under which conditions there is a stationary solution $((X_t, \sigma_t))$ to (3.1). In Section 3.2 we study the invertibility of this stationary solution. In Section 3.3 the function $h_t$, which was discussed in the Introduction of this paper, will be properly defined for model (3.1).

3.1. *Existence of a stationary solution.* Here and in the following section we suppress the parameter $\boldsymbol{\theta}$ in our notation because for the treatment of stationarity and invertibility $\boldsymbol{\theta}$ can be assumed as fixed. In order to discuss the stationarity issue, we first embed model (3.1) into an SRE. To this end, introduce ${}_k\boldsymbol{\sigma}_t^2 = (\sigma_t^2, \ldots, \sigma_{t-k+1}^2)^T$ and set ${}_k\mathbf{Z}_t = (Z_t, \ldots, Z_{t-k+1})^T$ for $k = \max(p, q)$. Then by substituting the $X_{t-i}$'s by $Z_{t-i}\sigma_{t-i}$ in the second equation of (3.1), we observe that $({}_k\boldsymbol{\sigma}_t^2)$ is a solution of the SRE

$$(3.2) \qquad \mathbf{s}_{t+1} = \psi_t(\mathbf{s}_t), \qquad t \in \mathbb{Z},$$

on $[0, \infty)^k$, where

$$(3.3) \qquad \psi_t(\mathbf{s}) = (g(\mathbf{s}^{1/2} \odot {}_k\mathbf{Z}_t, \mathbf{s}), s_1, \ldots, s_{k-1})^T,$$

with $\mathbf{s}^{1/2} = (s_1^{1/2}, \ldots, s_k^{1/2})^T$ for $\mathbf{s} \in [0, \infty)^k$; recall that $\odot$ stands for the Hadamard product, the componentwise multiplication of matrices or vectors of the same dimension. On the other hand, if a stationary sequence $(\mathbf{s}_t)$ is a solution of (3.2) and $\mathbf{s}_t$ is measurable with respect to $\mathcal{F}_{t-1}$, then the sequence $((\mathbf{s}_{t,1}^{1/2} Z_t, \mathbf{s}_{t,1}^{1/2}))$ (here $\mathbf{s}_{t,1}^{1/2}$ denotes the first coordinate of $\mathbf{s}_t^{1/2}$) is stationary and obeys the equations (3.1). Hence, for proving the existence and uniqueness of a stationary sequence $((X_t, \sigma_t))$ satisfying (3.1) and $\sigma_t^2$ being $\mathcal{F}_{t-1}$-measurable, we may focus on showing that there exists a unique stationary nonnegative solution $(\mathbf{s}_t)$ to (3.2) for which $\mathbf{s}_t$ is $\mathcal{F}_{t-1}$-measurable. After noticing that $(\psi_t)$ is stationary ergodic, an application of Theorem 2.8 with $(\phi_t) = (\psi_t)$ and $d$ the Euclidean metric results in the following proposition.

PROPOSITION 3.1. *Fix an arbitrary $\boldsymbol{\varsigma}_0^2 \in [0, \infty)^k$ and suppose that the following conditions hold for the stationary ergodic sequence $(\psi_t)$ of maps defined in* (3.3):

S.1 $\mathbb{E}[\log^+ |\psi_0(\boldsymbol{\varsigma}_0^2)|] < \infty$.

S.2 $\mathbb{E}[\log^+ \Lambda(\psi_0)] < \infty$ *and for some integer $r \geq 1$, $\mathbb{E}[\log \Lambda(\psi_0^{(r)})] < 0$.*

*Then the SRE* (3.2) *admits a unique stationary ergodic solution $({}_k\boldsymbol{\sigma}_t^2)$ such that ${}_k\boldsymbol{\sigma}_t^2$ is $\mathcal{F}_{t-1}$-measurable for every $t \in \mathbb{Z}$. The following almost sure representation is valid:*

$$(3.4) \qquad {}_k\boldsymbol{\sigma}_t^2 = \lim_{m \to \infty} \psi_{t-1} \circ \cdots \circ \psi_{t-m}(\boldsymbol{\varsigma}_0^2), \qquad t \in \mathbb{Z},$$

*where the latter limit is irrespective of $\boldsymbol{\varsigma}_0^2$. For any other solution $({}_k\tilde{\boldsymbol{\sigma}}_t^2)$ of the SRE* (3.2) *with index set $\mathbb{Z}$ or $\mathbb{N}$,*

$$(3.5) \qquad |{}_k\tilde{\boldsymbol{\sigma}}_t^2 - {}_k\boldsymbol{\sigma}_t^2| \xrightarrow{\text{e.a.s.}} 0, \qquad t \to \infty.$$

REMARK 3.2. In general, the stationary distribution induced by the SRE (3.2) is not known and, thus, perfect simulation of $((X_t, \sigma_t))_{t \in \mathbb{N}}$ is impossible. The following algorithm provides an approximation with an error decaying e.a.s.

(1) Take an initial value $\boldsymbol{\varsigma}_0^2 \in [0, \infty)^k$, set $_k\tilde{\boldsymbol{\sigma}}_0^2 = \boldsymbol{\varsigma}_0^2$, and generate $_k\tilde{\boldsymbol{\sigma}}_t^2$ according to (3.2).

(2) Set $(\tilde{X}_t, \tilde{\sigma}_t) = (_k\tilde{\boldsymbol{\sigma}}_{t,1}^2)^{1/2}(Z_t, 1)$, $t = 0, 1, \ldots$.

From Proposition 3.1, it is immediate that $|\tilde{\sigma}_t^2 - \sigma_t^2| \overset{\text{e.a.s.}}{\longrightarrow} 0$ as $t \to \infty$. Moreover, by virtue of Lemmas 2.1 and 2.3, we can conclude $|\tilde{X}_t - X_t| = |Z_t||\tilde{\sigma}_t - \sigma_t| \overset{\text{e.a.s.}}{\longrightarrow} 0$ as $t \to \infty$.

We now discuss two examples of conditionally heteroscedastic time series models of the form (3.1), which are well-known from the econometrics literature.

EXAMPLE 3.3. In its simplest form, the exponential GARCH model (EGARCH) of Nelson [27] has a volatility obeying the SRE

$$(3.6) \qquad \log \sigma_t^2 = \alpha + \beta \log \sigma_{t-1}^2 + \gamma Z_{t-1} + \delta |Z_{t-1}|, \qquad t \in \mathbb{Z},$$

where $\alpha, \gamma, \delta \in \mathbb{R}$ and $|\beta| < 1$. In other words, the sequence $(\log \sigma_t^2)$ constitutes a causal AR(1) process with innovations sequence $(\alpha + \gamma Z_{t-1} + \delta|Z_{t-1}|)$. Although EGARCH falls into the class of models defined by (3.1), we avoid writing the SRE (3.6) in the form $\sigma_t^2 = g(X_{t-1}, \sigma_{t-1}^2)$ since this would unnecessarily complicate the question of stationarity. Instead, we apply Theorem 2.8 to the SRE

$$(3.7) \qquad \log \sigma_{t+1}^2 = \psi_t(\log \sigma_t^2), \qquad t \in \mathbb{Z},$$

where

$$\psi_t(s) = \alpha + \beta s + \gamma Z_t + \delta |Z_t|, \qquad t \in \mathbb{Z}.$$

We recognize that $|\beta| < 1$ together with $\mathbb{E}[\log^+(\alpha + \gamma Z_0 + \delta|Z_0|)] < \infty$ is a sufficient condition for the existence of a unique stationary solution to (3.7); since $\mathbb{E}Z_0^2 = 1$, the innovations $(\alpha + \gamma Z_{t-1} + \delta|Z_{t-1}|)$ automatically have a finite positive logarithmic moment. By taking the limit of the backward iterates associated with the SRE (3.7), it is straightforward to see that the unique stationary solution $(\log \sigma_t^2)$ of (3.7) has the almost sure representation

$$(3.8) \quad \log \sigma_t^2 = \alpha(1-\beta)^{-1} + \sum_{k=0}^{\infty} \beta^k (\gamma Z_{t-1-k} + \delta|Z_{t-1-k}|), \qquad t \in \mathbb{Z},$$

which is well-known from the theory of the causal AR(1) process; see, for example, [11].

EXAMPLE 3.4. In the AGARCH$(p,q)$ model [asymmetric GARCH$(p,q)$], which was independently proposed by Ding, Granger and Engle [14] and Zakoïan [33], the squared volatility is of form

$$(3.9) \qquad \sigma_t^2 = \alpha_0 + \sum_{i=1}^{p} \alpha_i (|X_{t-i}| - \gamma X_{t-i})^2 + \sum_{j=1}^{q} \beta_j \sigma_{t-j}^2, \qquad t \in \mathbb{Z},$$

where $\alpha_0 > 0$, $\alpha_i, \beta_j \geq 0$ and $|\gamma| \leq 1$. Note that $\gamma = 0$ in AGARCH$(p,q)$ corresponds to a GARCH$(p,q)$ process. One may apply Proposition 3.1 in order to study the stationarity of AGARCH, but the resulting sufficient conditions will not be very enlightening. Instead, we adapt the particular methods developed in [6] for the treatment of GARCH$(p,q)$ and derive an SRE with i.i.d. linear random transition maps. This leads to a more insightful theory. To this end, consider the associated sequence $(\mathbf{Y}_t)$ given by

$$(3.10) \qquad \begin{aligned} \mathbf{Y}_t = (&\sigma_t^2, \ldots, \sigma_{t-q+1}^2, \\ &(|X_{t-1}| - \gamma X_{t-1})^2, \ldots, (|X_{t-p+1}| - \gamma X_{t-p+1})^2)^T. \end{aligned}$$

One easily verifies that $X_t = \sigma_t Z_t$ together with (3.9) implies

$$(3.11) \qquad \mathbf{Y}_{t+1} = \mathbf{A}_t \mathbf{Y}_t + \mathbf{B}_t, \qquad t \in \mathbb{Z},$$

where the $(p+q-1) \times (p+q-1)$-matrix valued i.i.d. sequence $(\mathbf{A}_t)$ and the vectors $(\mathbf{B}_t)$ are defined through

$$(3.12) \quad \mathbf{A}_t = \begin{pmatrix} \alpha_1(|Z_t| - \gamma Z_t)^2 + \beta_1 & \beta_2 & \cdots & \beta_{q-1} & \beta_q & \alpha_2 & \alpha_3 & \cdots & \alpha_{p-1} & \alpha_p \\ 1 & 0 & \cdots & 0 & 0 & 0 & 0 & \cdots & 0 & 0 \\ 0 & 1 & \cdots & 0 & 0 & 0 & 0 & \cdots & 0 & 0 \\ \vdots & & \ddots & \vdots & \vdots & \vdots & \vdots & \ddots & \vdots & \vdots \\ 0 & 0 & \cdots & 1 & 0 & 0 & 0 & \cdots & 0 & 0 \\ (|Z_t| - \gamma Z_t)^2 & 0 & \cdots & 0 & 0 & 0 & 0 & \cdots & 0 & 0 \\ 0 & 0 & \cdots & 0 & 0 & 1 & 0 & \cdots & 0 & 0 \\ \vdots & & \vdots & \vdots & \vdots & \vdots & \ddots & \ddots & \vdots & \vdots \\ 0 & 0 & \cdots & 0 & 0 & 0 & 0 & \cdots & 1 & 0 \end{pmatrix},$$

$$\mathbf{B}_t = (\alpha_0, 0, \ldots, 0)^T \in \mathbb{R}^{p+q-1}.$$

A straightforward argument yields that there is a unique stationary AGARCH$(p,q)$ process $((X_t, \sigma_t))$ if and only if the SRE (3.11) has a unique stationary nonnegative solution such that $\mathbf{Y}_t$ is $\mathcal{F}_{t-1}$-measurable. By employing the identical arguments as in the proof of Theorem 1.3 in [6], the SRE (3.11) has a unique stationary solution if and only if the top Lyapunov exponent associated with $(\mathbf{A}_t)$ is strictly negative, that is,

$$(3.13) \qquad \rho = \lim_{n \to \infty} \frac{1}{n+1} \mathbb{E}(\log \|\mathbf{A}_0 \cdots \mathbf{A}_{-n+1}\|_{\mathrm{op}}) < 0.$$

Here $\|\cdot\|_{\text{op}}$ denotes the matrix operator norm with respect to the Euclidean norm on $\mathbb{R}^{p+q-1}$. Moreover, this stationary solution is ergodic and measurable with respect to $\mathcal{F}_{t-1}$ with almost sure representation

$$(3.14) \qquad \mathbf{Y}_t = \sum_{k=0}^{\infty} \left( \prod_{l=1}^{k} \mathbf{A}_{t-l} \right) \mathbf{B}_{t-k}, \qquad t \in \mathbb{Z}.$$

Here by convention $\prod_{l=1}^{0} \cdot \equiv 1$. Furthermore, one demonstrates along the lines of the proof of Corollary 2.2 in [6] and the remark following it that

$$(3.15) \qquad \sum_{i=1}^{p} \alpha_i \mathbb{E}[(|Z_0| - \gamma Z_0)^2] + \sum_{j=1}^{q} \beta_j < 1$$

is sufficient for stationarity and implies $\mathbb{E}X_0^2 < \infty$. This property generalizes the necessary and sufficient condition of [4] for weak stationarity in GARCH$(p,q)$ processes to the AGARCH class. Arguing as in Corollary 2.3 of [6], one finds that

$$(3.16) \qquad \sum_{j=1}^{q} \beta_j < 1$$

is necessary for stationarity of AGARCH$(p,q)$. To the best of our knowledge, these results do not appear in the literature, for which reason we summarize them in the form of a theorem.

THEOREM 3.5. *The AGARCH$(p,q)$ equations $X_t = \sigma_t Z_t$ and (3.9) admit a unique stationary solution if and only if $(\mathbf{A}_t)$ has top Lyapunov exponent $\rho < 0$. Moreover, this stationary solution $((X_t, \sigma_t))$ is ergodic, $\sigma_t$ is $\mathcal{F}_{t-1}$-measurable, and the vector $\mathbf{Y}_t$ given in (3.10) has almost sure representation (3.14). The restriction (3.16) on the parameters is necessary for stationarity, and the condition (3.15) is equivalent to the existence and uniqueness of a stationary AGARCH$(p,q)$ process with finite second moment marginal distribution.*

REMARK 3.6. The sufficiency of $\rho < 0$ for the existence of a stationary solution to (3.11) is a direct consequence of Theorem 2.8. Indeed, the SRE (3.11) can be written as $\mathbf{Y}_{t+1} = \psi_t(\mathbf{Y}_t)$, where $\psi_t(\mathbf{y}) = \mathbf{A}_t \mathbf{y} + \mathbf{B}_t$. It is clear that $\Lambda(\psi_0^{(m)}) \leq \|\mathbf{A}_0^{(m)}\|_{\text{op}} = \|\mathbf{A}_0 \cdots \mathbf{A}_{-m+1}\|_{\text{op}}$ for all $m \in \mathbb{N}$. By definition of the top Lyapunov exponent, there is an $r \geq 1$ such that $\mathbb{E}(\log \|\mathbf{A}_0^{(r)}\|_{\text{op}}) < 0$ and, hence, the conditions of Theorem 2.8 are met with $(\phi_t) = (\psi_t)$ and $d$ the Euclidean norm; note that the technical assumption S.1 and $\mathbb{E}[\log^+ \Lambda(\psi_0)] = \mathbb{E}(\log^+ \|\mathbf{A}_0\|_{\text{op}}) < \infty$ trivially follow from $\mathbb{E}\|\mathbf{A}_0\|_{\text{op}} < \infty$, which is a consequence of $\mathbb{E}Z_0^2 = 1$ and the fact that every matrix norm is equivalent to the

Frobenius norm. The definition (3.13) of the top Lyapunov exponent and $\rho < 0$ imply

$$\|\mathbf{A}_{t-1} \cdots \mathbf{A}_{t-m}\|_{\mathrm{op}} \xrightarrow{\text{e.a.s.}} 0, \qquad m \to \infty.$$

The latter property can be used for determining the limit of the backward iterates associated with (3.11), resulting in the representation (3.14). We also mention that $\mathbb{E}|\mathbf{Y}_0|^{\tilde{q}} < \infty$ for some $\tilde{q} > 0$. Indeed, $\mathbb{E}(\log \|\mathbf{A}_0^{(r)}\|_{\mathrm{op}}) < 0$ together with Lemma 2.4 implies the existence of $\tilde{q} > 0$ and $0 < \eta < 1$ such that

$$(3.17) \qquad \mathbb{E}\left( \prod_{\ell=1}^{m} \|\mathbf{A}_{-(\ell-1)r}^{(r)}\|_{\mathrm{op}} \right)^{\tilde{q}} = \eta^m, \qquad m \geq 1.$$

Without loss of generality, $\tilde{q} \leq 1$. Since $\mathbb{E}\|\mathbf{A}_0\|_{\mathrm{op}} < \infty$, also $\mathbb{E}\|\mathbf{A}_0\|_{\mathrm{op}}^{\tilde{q}} < \infty$. The identity (3.17) together with the facts that $\|\cdot\|_{\mathrm{op}}$ is submultiplicative and $(\mathbf{A}_t)$ is i.i.d. demonstrates

$$\mathbb{E}\|\mathbf{A}_0^{(k)}\|_{\mathrm{op}}^{\tilde{q}} \leq (\mathbb{E}\|\mathbf{A}_0^{(r)}\|_{\mathrm{op}}^{\tilde{q}})^{[k/r]} (\mathbb{E}\|\mathbf{A}_0\|_{\mathrm{op}}^{\tilde{q}})^{r-r[k/r]} \leq c\eta^{[k/r]}, \qquad k \geq 0,$$

where $c = \max(1, (\mathbb{E}\|\mathbf{A}_0\|_{\mathrm{op}}^{\tilde{q}})^r)$. Finally, by an application of the Minkowski inequality,

$$(3.18) \qquad \mathbb{E}\|\mathbf{Y}_1\|_{\mathrm{op}}^{\tilde{q}} \leq \sum_{k=0}^{\infty} \mathbb{E}\|\mathbf{A}_0^{(k)}\|_{\mathrm{op}}^{\tilde{q}} \mathbb{E}|\mathbf{B}_0|^{\tilde{q}} \leq cc_0^{\tilde{q}} \sum_{k=0}^{\infty} \eta^{[k/r]} < \infty.$$

### 3.2. *Invertibility.*

For real-life data sets which we assume to be generated by a model of type (3.1), the volatility $\sigma_t$ will be unobservable. In such a case, it is natural to approximate the unobservable squared volatilities from data $X_{-p+1}, X_{-p+2}, \ldots$ in the following way:

*Initialization*:

(1) Set $\hat{\boldsymbol{\sigma}}_0^2 = \boldsymbol{\varsigma}_0^2$, where $\boldsymbol{\varsigma}_0^2 \in [0, \infty)^q$ is an arbitrarily chosen vector.

*Recursion step*:

(2) Let $\hat{\boldsymbol{\sigma}}_{t+1}^2 = \phi_t(\hat{\boldsymbol{\sigma}}_t^2)$ for $t = 0, 1, 2, \ldots$, where the random functions $\phi_t$ are defined as

$$(3.19) \quad \phi_t(\mathbf{s}) = (g(\mathbf{X}_t, \mathbf{s}), s_1, \ldots, s_{q-1})^T, \qquad \mathbf{s} = (s_1, \ldots, s_q)^T \in [0, \infty)^q.$$

In the context of our nonlinear model we say that the unique stationary ergodic solution $((X_t, \sigma_t))$ to the equations (3.1) is *invertible* [or model (3.1) is invertible] if

$$|\hat{\boldsymbol{\sigma}}_t^2 - \boldsymbol{\sigma}_t^2| \xrightarrow{\mathbb{P}} 0, \qquad t \to \infty.$$

In other words, invertibility guarantees that the above algorithm converges. There is a second interpretation, which at the same time clarifies the relationship between invertibility in ARMA models and our notion. Note that $(\hat{\boldsymbol{\sigma}}_t^2)_{t \in \mathbb{N}}$ is a solution of the SRE

$$\mathbf{s}_{t+1} = \phi_t(\mathbf{s}_t), \qquad t \in \mathbb{N},$$

on $[0, \infty)^q$ and recall that the backward iterates associated with $(\phi_t)$ are defined by

$$\boldsymbol{\sigma}_{t,0}^2 = \boldsymbol{\varsigma}_0^2,$$
$$\boldsymbol{\sigma}_{t,m}^2 = \phi_{t-1} \circ \cdots \circ \phi_{t-m}(\boldsymbol{\varsigma}_0^2), \qquad m \geq 1,$$

for $t \in \mathbb{Z}$. Observe the relationship $\hat{\boldsymbol{\sigma}}_t^2 = \boldsymbol{\sigma}_{t,t}^2$ for all $t \geq 0$. Furthermore, since we suppose that $((X_t, \sigma_t))$ is stationary, $((\boldsymbol{\sigma}_{t,m}^2, \boldsymbol{\sigma}_t^2))$ is stationary for every fixed $m \geq 0$ (Proposition 2.5). Thus,

$$\boldsymbol{\sigma}_{t,m}^2 - \boldsymbol{\sigma}_t^2 \stackrel{\mathrm{d}}{=} \boldsymbol{\sigma}_{m,m}^2 - \boldsymbol{\sigma}_m^2 = \hat{\boldsymbol{\sigma}}_m^2 - \boldsymbol{\sigma}_m^2, \qquad m \geq 0.$$

Therefore, invertibility is equivalent to

$$(3.20) \qquad \boldsymbol{\sigma}_{t,m}^2 \stackrel{\mathbb{P}}{\longrightarrow} \boldsymbol{\sigma}_t^2, \qquad m \to \infty,$$

for every $t \in \mathbb{Z}$. Relation (3.20) together with Proposition 2.6 and a subsequence argument implies the existence of a measurable function $f$ such that

$$\boldsymbol{\sigma}_t^2 = f(X_{t-1}, X_{t-2}, \ldots) \qquad \text{a.s.}$$

for every $t \in \mathbb{Z}$. If we impose $\sigma_t^2 > 0$ a.s. for the model (3.1), the invertibility allows us to represent $Z_t$ as a function of the past and present observations $\{X_{t-k} | k \geq 0\}$. Compare this with the notion of invertibility in ARMA models. Recall that an ARMA$(p, q)$ model with parameters $(\alpha_1, \ldots, \alpha_p, \beta_1, \ldots, \beta_q)^T$ is a stochastic process $(X_t)$ which obeys the difference equation

$$X_t = \sum_{i=1}^p \alpha_i X_{t-i} + \sum_{j=1}^q \beta_j Z_{t-j} + Z_t, \qquad t \in \mathbb{Z},$$

where $(Z_t)$ is a given white noise sequence with mean zero and finite variance $\sigma^2$; see [11]. Invertibility is defined as follows: there exists an absolutely summable sequence of constants $(\pi_j)_{j \geq 0}$ such that $Z_t = \sum_{j=0}^{\infty} \pi_j X_{t-j}$ a.s., or, in other words, the innovation at time $t$ is a linear functional of the past and present observations $\{X_{t-k} | k \geq 0\}$. If the two characteristic polynomials $\alpha(z) = 1 - \sum_{i=1}^p \alpha_i z^i$ and $\beta(z) = 1 + \sum_{j=1}^q \beta_j z^j$ have no common zeros, then the ARMA process is invertible if and only if $\beta(z) \neq 0$ for all $z \in \mathbb{C}$ such that $|z| \leq 1$; see, for example, Theorem 3.1.2 in [11]. For nonlinear time series models, however, the invertibility issue can be a hard problem. Our notion of invertibility is an adaptation of the notion introduced by Granger and Andersen [16] in the context of a general nonlinear autoregressive moving average model. The following proposition is an immediate consequence of Theorem 2.8 with $d$ the Euclidean metric.

PROPOSITION 3.7. *Assume that there exists a unique stationary ergodic solution $((X_t, \sigma_t))$ to the model* (3.1). *Suppose, in addition, for $(\phi_t)$ given in* (3.19):

S.1 $\mathbb{E}[\log^+ |\phi_0(\varsigma_0^2)|] < \infty$.

S.2 $\mathbb{E}[\log^+ \Lambda(\phi_0)] < \infty$ *and for some integer $r \geq 1$, $\mathbb{E}[\log \Lambda(\phi_0^{(r)})] < 0$.*

*If there exists $r \geq 1$ with $\mathbb{E}[\log \Lambda(\phi_0^{(r)})] < 0$, then $((X_t, \sigma_t))$ is invertible. In particular, irrespective of $\hat{\varsigma}_0^2$,*

$$|\hat{\boldsymbol{\sigma}}_t^2 - \boldsymbol{\sigma}_t^2| \xrightarrow{\text{e.a.s.}} 0, \qquad t \to \infty.$$

EXAMPLE 3.8 (Continuation of Example 3.3). Assume $0 \leq \beta < 1$, $\gamma + \delta \geq 0$, $\delta - \gamma \geq 0$ for the parameters of EGARCH. This restriction also seems reasonable from an economics point of view. One expects a positive relationship between volatilities on successive days, that is, $\beta \geq 0$. The squared volatility $\sigma_t^2$ as a function of $Z_{t-1}$ should be nondecreasing on the positive real line (i.e., $\gamma + \delta \geq 0$) and nonincreasing on the negative real line (i.e., $\delta - \gamma \geq 0$). Altogether $\gamma z + \delta|z| \geq 0$ for all $z \in \mathbb{R}$. As we mentioned in Example 3.3, it is beneficial to consider the SRE for $(\log \sigma_t^2)$, which slightly differs from the general setup of this paper. From (3.8), one concludes that $\log \sigma_t^2 \geq \alpha(1 - \beta)^{-1}$ so that we may interpret (3.6) as an SRE

$$(3.21) \qquad \log \sigma_{t+1}^2 = \phi_t(\log \sigma_t^2), \qquad t \in \mathbb{Z},$$

on the (restricted) set $I = [\alpha(1-\beta)^{-1}, \infty)$, where

$$\phi_t(s) = \alpha + \beta s + (\gamma X_t + \delta|X_t|) \exp(-s/2).$$

Let us determine $\Lambda(\phi_0)$. Since $\phi_0$ is continuously differentiable, $\Lambda(\phi_0) = \sup_{s \in I} |\phi_0'(s)|$. Maximizing

$$|\phi_0'(s)| = |\beta - 2^{-1}(\gamma X_0 + \delta|X_0|) \exp(-s/2)|$$

over $I$, we obtain

$$(3.22) \quad \Lambda(\phi_0) = \max(\beta, 2^{-1} \exp(-\alpha(1-\beta)^{-1}/2)(\gamma X_0 + \delta|X_0|)) - \beta).$$

It does not seem to be possible to derive a tractable expression for $\Lambda(\phi_0^{(r)})$ when $r > 1$ and it is not clear how to find a bound sharper than the trivial bound $\Lambda(\phi_0^{(r)}) \leq \Lambda(\phi_0) \cdots \Lambda(\phi_{-r+1})$. Recalling the representation (3.8) and substituting $X_0$ by $\sigma_0 Z_0$ in (3.22), the condition $\mathbb{E}[\log \Lambda(\phi_0)] < 0$, which implies invertibility, reads

$$
\text{(3.23)} \quad
\begin{aligned}
\mathbb{E}\bigg[ \log \max \bigg\{ & \beta, 2^{-1} \exp \bigg( 2^{-1} \sum_{k=0}^{\infty} \beta^k (\gamma Z_{-k-1} + \delta |Z_{-k-1}|) \bigg) \\
& \times (\gamma Z_0 + \delta |Z_0|) - \beta \bigg\} \bigg] < 0.
\end{aligned}
$$

In practice, one would have to rely on simulation methods for its verification.

PROPOSITION 3.9. *Let $((X_t, \sigma_t))$ be a stationary ergodic of the EGARCH process with the parameter constraints $0 \le \beta < 1$, $\delta \ge |\gamma|$. Then the condition (3.23) is sufficient for invertibility since*

$$
\text{(3.24)} \qquad |\hat{\sigma}_t^2 - \sigma_t^2| \overset{\text{e.a.s.}}{\longrightarrow} 0, \qquad t \to \infty.
$$

PROOF. By an application of Theorem 2.8 to the SRE (3.21), one has $|\log \hat{\sigma}_t^2 - \log \sigma_t^2| \overset{\text{e.a.s.}}{\longrightarrow} 0$. Taking the expectation in (3.8) and accounting for $\mathbb{E}|Z_t| \le (\mathbb{E}Z_t^2)^{1/2} = 1$, we obtain $\mathbb{E}(\log \sigma_0^2) \le (\alpha + \delta)(1 - \beta)^{-1} < \infty$. Consequently, (3.24) is valid by virtue of Lemma 2.3. □

REMARK 3.10. The case $\beta = 0$ leads to a simpler condition and illustrates that there exist invertible EGARCH models. For $\beta = 0$, the condition (3.23) becomes $-\log 2 + (\delta/2)\mathbb{E}|Z_0| + \mathbb{E}[\log(\gamma Z_0 + \delta |Z_0|)] < 0$. Note that $\delta \le 1$ implies the latter condition. Indeed, using $\mathbb{E}|Z_0| \le (\mathbb{E}Z_0^2)^{1/2} = 1$ and $\mathbb{E}(\gamma Z_0 + \delta |Z_0|) \le \delta \le 1$, Jensen's inequality yields $\mathbb{E}[\log(\gamma Z_0 + \delta |Z_0|)] \le 0$. From the relation $1/2 < \log 2$, one obtains $(\delta/2)\mathbb{E}|Z_0| - \log 2 < 0$, which proves the assertion.

EXAMPLE 3.11 (Continuation of Example 3.4). For the AGARCH$(p, q)$ model, the random maps (3.19) are affine linear, that is,

$$
\phi_t(\mathbf{s}) = \bigg( \alpha_0 + \sum_{i=1}^{p} \alpha_i (|X_{t+1-i}| - \gamma X_{t+1-i})^2 \bigg) \mathbf{e}_1 + \mathbf{C}\mathbf{s}, \qquad \mathbf{s} \in [0, \infty)^q,
$$

where

$$
\mathbf{e}_1 = (1, 0, \ldots, 0)^T \in \mathbb{R}^q,
$$

$$
\mathbf{C} = \begin{pmatrix}
\beta_1 & \beta_2 & \cdots & \cdots & \beta_{q-1} & \beta_q \\
1 & 0 & \cdots & \cdots & 0 & 0 \\
0 & 1 & 0 & \cdots & \cdots & 0 \\
\vdots & \ddots & \ddots & \ddots & \vdots & \vdots \\
0 & \cdots & 0 & 1 & 0 & 0 \\
0 & \cdots & \cdots & 0 & 1 & 0
\end{pmatrix} \in \mathbb{R}^{q \times q}.
$$

Observe that $\Lambda(\phi_0^{(r)}) \leq \|\mathbf{C}^r\|_{\mathrm{op}}$ for any $r \in \mathbb{N}$ (the inequality in the latter relation is a consequence of the fact that the domain of $\phi_t$ is but a subset of $\mathbb{R}^q$). In order to prove $\log \|\mathbf{C}^r\|_{\mathrm{op}} < 0$ for large enough $r$, first recall that a necessary condition for stationarity in $\mathrm{AGARCH}(p,q)$ is $\beta_\Sigma := \sum_{j=1}^q \beta_j < 1$. Arguing by recursion on $q$ and expanding the determinant with respect to the last column, it can be shown that $\mathbf{C}$ has characteristic polynomial $p(\lambda) = \det(\lambda \mathbf{I} - \mathbf{C}) = \lambda^q(1 - \sum_{j=1}^q \beta_j \lambda^{-i})$. If $|\lambda| > \beta_\Sigma^{1/q}$, then by repeated application of the triangle inequality together with $0 < \beta_\Sigma < 1$,

$$|p(\lambda)| \geq 1 - \sum_{j=1}^q \beta_j |\lambda|^{-j} > 1 - \sum_{j=1}^q \beta_j \beta_\Sigma^{-j/q} \geq 1 - \beta_\Sigma^{-1} \sum_{j=1}^q \beta_j = 0.$$

Consequently, the matrix $\mathbf{C}$ has spectral radius strictly smaller than $\beta_\Sigma^{1/q}$. By the Jordan normal decomposition, this entails

$$(3.25) \qquad \|\mathbf{C}^r\|_{\mathrm{op}} \leq c\beta_\Sigma^{r/q}, \qquad r \geq 0,$$

where $c$ is a constant depending on $q$ and $\beta_\Sigma$, and thus, $\log \Lambda(\phi_0^{(r)}) \to -\infty$. This means that the conditions of Proposition 3.7 are met for stationary $\mathrm{AGARCH}(p,q)$ processes. Consequently, *every* stationary $\mathrm{AGARCH}(p,q)$ process is automatically invertible. It is also possible to give an explicit representation of $\sigma_t^2$ in terms of $(X_{t-1}, X_{t-2}, \ldots)$; see equation (5.6) below.

3.3. *Definition of the function* $\mathbf{h}_t$. In this section we do not suppress the parameter $\boldsymbol{\theta}$ any more, that is, we write $\sigma_{t+1}^2 = g_{\boldsymbol{\theta}}(\mathbf{X}_t, \boldsymbol{\sigma}_t^2)$. Assume that $\boldsymbol{\theta}$ belongs to some compact parameter space $K \subset \mathbb{R}^d$ and that $((X_t, \sigma_t))$ is the unique stationary ergodic solution to the model (3.1) with unknown true parameter $\boldsymbol{\theta} = \boldsymbol{\theta}_0$. For any initial value $\boldsymbol{\varsigma}_0^2 \in [0, \infty)^q$, we define the following random vector functions $\hat{\mathbf{h}}_t$ on $K$:

$$(3.26) \qquad \hat{\mathbf{h}}_t = \begin{cases} \boldsymbol{\varsigma}_0^2, & t = 0, \\ \Phi_{t-1}(\hat{\mathbf{h}}_{t-1}), & t \geq 1, \end{cases}$$

where the random maps $\Phi_t : \mathbb{C}(K, [0,\infty)^q) \to \mathbb{C}(K, [0,\infty)^q)$ are given by

$$[\Phi_t(\mathbf{s})](\boldsymbol{\theta}) = (g_{\boldsymbol{\theta}}(\mathbf{X}_t, \mathbf{s}(\boldsymbol{\theta})), s_1(\boldsymbol{\theta}), \ldots, s_{q-1}(\boldsymbol{\theta}))^T, \qquad t \in \mathbb{Z}.$$

We can regard $\hat{\mathbf{h}}_t(\boldsymbol{\theta}) = (\hat{h}_t(\boldsymbol{\theta}), \ldots, \hat{h}_{t-q+1}(\boldsymbol{\theta}))^T$ as an "estimate" of the squared volatility vector $\boldsymbol{\sigma}_t^2$ under the parameter hypothesis $\boldsymbol{\theta}$, which is based on the data $X_{-p+1}, \ldots, X_t$. Also observe that $\hat{\mathbf{h}}_t(\boldsymbol{\theta}_0) = \hat{\boldsymbol{\sigma}}_t^2$ for all $t \in \mathbb{N}$.

For establishing the consistency of the QMLE, it is essential that one can approximate $(\hat{h}_t)_{t \in \mathbb{N}}$ by a stationary ergodic sequence $(h_t)_{t \in \mathbb{N}}$ such that the error $\hat{h}_t - h_t$ converges to zero sufficiently fast as $t \to \infty$ and such that $(h_t(\boldsymbol{\theta})) = (\sigma_t^2)$ a.s. if and only if $\boldsymbol{\theta} = \boldsymbol{\theta}_0$; see Theorem 4.1. In particular, these requirements on $(h_t)_{t \in \mathbb{N}}$ comprehend the invertibility of the time series $(X_t)$ because $\hat{h}_t(\boldsymbol{\theta}_0) = \hat{\sigma}_t^2$. Invertibility is a common assumption in the classical theory of parameter estimation in ARMA time series. Only recently have certain aspects of estimation in noninvertible linear time series been studied; see, for example, [13] or [10].

For finding a candidate for such an approaching sequence, first observe that $(\hat{\mathbf{h}}_t)_{t \in \mathbb{N}}$ is a solution of the SRE

$$(3.27) \qquad \mathbf{s}_{t+1} = \Phi_t(\mathbf{s}_t), \qquad t \in \mathbb{N},$$

on $\mathbb{C}(K, [0, \infty)^q)$ (provided certain regularity assumptions on $g_{\boldsymbol{\theta}}$ are fulfilled). Theorem 2.8 at hand, it is clear that if $\Phi_0$ or one of its $r$-fold convolutions is a contraction on average, then the unique stationary ergodic solution of the SRE (3.27) with index set $\mathbb{Z}$ provides the desired sequence $(\mathbf{h}_t)$. We summarize our findings in a proposition.

PROPOSITION 3.12. *Assume that model* (3.1) *admits a unique stationary ergodic solution* $((X_t, \sigma_t))$ *and that the map* $(\boldsymbol{\theta}, \mathbf{s}) \mapsto g_{\boldsymbol{\theta}}(\mathbf{x}, \mathbf{s})$ *is continuous for every* $\mathbf{x} \in \mathbb{R}^p$, *which implies that* $(\Phi_t)$ *is a stationary ergodic sequence of mappings* $\mathbb{C}(K, [0, \infty)^q) \to \mathbb{C}(K, [0, \infty)^q)$. *We suppose the following conditions hold:*

1. $\mathbb{E}(\log^+ \|\Phi_0(\boldsymbol{\varsigma}_0^2)\|_K) < \infty$.
2. $\mathbb{E}[\log^+ \Lambda(\Phi_0)] < \infty$ *and there exists an integer* $r \geq 1$ *such that* $\mathbb{E}[\log \Lambda(\Phi_0^{(r)})] < 0$.

*Then the SRE* (3.27) *(with the index set* $\mathbb{N}$ *replaced by* $\mathbb{Z}$*) has a unique stationary solution* $(\mathbf{h}_t)$, *which is ergodic. For every* $t \in \mathbb{Z}$, *the random elements* $\mathbf{h}_t$ *are* $\mathcal{F}_{t-1}$*-measurable and* $\mathbf{h}_t(\boldsymbol{\theta}_0) = \boldsymbol{\sigma}_t^2$ *a.s. Moreover,*

$$(3.28) \qquad \|\hat{\mathbf{h}}_t - \mathbf{h}_t\|_K \xrightarrow{\text{e.a.s.}} 0, \qquad t \to \infty.$$

PROOF. The existence and uniqueness of a stationary ergodic solution to (3.27) and the limit relation (3.28) are a direct consequence of Theorem 2.8. It remains to give arguments for the remaining assertions. The backward iterates associated with (3.27) are given by

$$(3.29) \qquad \mathbf{h}_{t,m} = \begin{cases} \boldsymbol{\varsigma}_0^2, & m = 0, \\ \Phi_{t-1} \circ \cdots \circ \Phi_{t-m}(\boldsymbol{\varsigma}_0^2), & m \geq 1. \end{cases}$$

This reveals that they are of the form $\mathbf{h}_{t,m} = f_m(X_{t-1}, X_{t-2}, \ldots)$ for certain measurable maps $f_m$. Since $\mathbf{h}_{t,m} \xrightarrow{\text{a.s.}} \mathbf{h}_t$ as $t \to \infty$, an application of Proposition 2.6 shows that $\mathbf{h}_t = f(X_{t-1}, X_{t-2}, \ldots)$ a.s., where $f$ is measurable. From this we conclude that $\mathbf{h}_t$ is $\mathcal{F}_{t-1}$-measurable for every $t$. The relation $\mathbf{h}_t(\boldsymbol{\theta}_0) = \boldsymbol{\sigma}_t^2$ a.s. follows from

$$\mathbf{h}_{t,m}(\boldsymbol{\theta}_0) - \boldsymbol{\sigma}_t^2 = \boldsymbol{\sigma}_{t,m}^2 - \boldsymbol{\sigma}_t^2 \stackrel{\text{d}}{=} \hat{\boldsymbol{\sigma}}_m^2 - \boldsymbol{\sigma}_m^2,$$

together with $|\hat{\boldsymbol{\sigma}}_m^2 - \boldsymbol{\sigma}_m^2| \overset{\text{e.a.s.}}{\longrightarrow} 0$, as shown in Proposition 3.7.    □

**4. Consistency of the QMLE.** Suppose we observe data $X_{-p+1}, \ldots, X_0$, $X_1, \ldots, X_n$ generated by the model (3.1) with $\boldsymbol{\theta}_0$ as the true parameter; here we have to emphasize that by a shift of the index we can always assume that the data $X_{-p+1}, \ldots, X_0$ are available to us. By this convention, $\hat{h}_1$ is then well defined. As in the GARCH(1, 1) case discussed in the Introduction we define the conditional Gaussian likelihood by

$$(4.1) \qquad \hat{L}_n(\boldsymbol{\theta}) = -\frac{1}{2} \sum_{t=1}^n \left( \frac{X_t^2}{\hat{h}_t(\boldsymbol{\theta})} + \log \hat{h}_t(\boldsymbol{\theta}) \right),$$

where $\hat{h}_t$ is the first coordinate of the random vector function $\hat{\mathbf{h}}_t = (\hat{h}_t, \ldots, \hat{h}_{t-q+1})^T$ defined in (3.26). The function $\hat{\mathbf{h}}_t(\boldsymbol{\theta})$ serves as an estimate of the vector $\boldsymbol{\sigma}_t^2 = (\sigma_t^2, \ldots, \sigma_{t-q+1}^2)$ under the parameter hypothesis $\boldsymbol{\theta}$; see Section 3.3. The QMLE $\hat{\boldsymbol{\theta}}_n$ maximizes $\hat{L}_n$ on $K$, where $K$ is an appropriately chosen compact subset of the parameter space $\Theta$, that is,

$$(4.2) \qquad \hat{\boldsymbol{\theta}}_n = \underset{\boldsymbol{\theta} \in K}{\arg\max} \, \hat{L}_n(\boldsymbol{\theta}).$$

Assume that the conditions of Proposition 3.12 are satisfied. Then $\|\hat{\mathbf{h}}_t - \mathbf{h}_t\|_K \overset{\text{e.a.s.}}{\longrightarrow} 0$ as $t \to \infty$, where $(\mathbf{h}_t)$ is the unique stationary solution of the SRE $\mathbf{s}_{t+1} = \Phi_t(\mathbf{s}_t)$, $t \in \mathbb{Z}$. Then, following the ideas presented in the Introduction, we define

$$(4.3) \qquad L_n(\boldsymbol{\theta}) = -\frac{1}{2} \sum_{t=1}^n \left( \frac{X_t^2}{h_t(\boldsymbol{\theta})} + \log h_t(\boldsymbol{\theta}) \right)$$

together with

$$(4.4) \qquad \tilde{\boldsymbol{\theta}}_n = \underset{\boldsymbol{\theta} \in K}{\arg\max} \, L_n(\boldsymbol{\theta}).$$

From a theoretical point of view, it is more convenient to work with $(L_n)$ because $(X_t^2/h_t + \log h_t)$ is stationary ergodic, whereas $(X_t^2/\log \hat{h}_t + \log \hat{h}_t)_{t \in \mathbb{N}}$ is not. In what follows, we give a set of conditions which implies the strong consistency of $\hat{\boldsymbol{\theta}}_n$:

C.1 Model (3.1) with $\boldsymbol{\theta} = \boldsymbol{\theta}_0$ admits a unique stationary ergodic solution $((X_t, \sigma_t))$ with $\mathbb{E}(\log^+ \sigma_0^2) < \infty$.

C.2 The conditions of Proposition 3.12 are satisfied for a compact set $K \subset \Theta$ with $\boldsymbol{\theta}_0 \in K$.

C.3 The class of functions $\{g_{\boldsymbol{\theta}} | \boldsymbol{\theta} \in K\}$ is uniformly bounded from below, that is, there exists a constant $\underline{g} > 0$ such that $g_{\boldsymbol{\theta}}(\mathbf{x}, \mathbf{s}) \geq \underline{g}$ for all $(\mathbf{x}, \mathbf{s}) \in \mathbb{R}^p \times [0, \infty)^q$ and $\boldsymbol{\theta} \in K$.

C.4  The following identifiability condition holds on $K$: for all $\boldsymbol{\theta} \in K$,

$$h_0(\boldsymbol{\theta}) \equiv \sigma_0^2 \text{ a.s.} \quad \text{if and only if} \quad \boldsymbol{\theta} = \boldsymbol{\theta}_0.$$

C.5  The random elements $\sigma_0^2/h_0$ and $\log h_0$ have a finite expected norm,

$$\mathbb{E}\left\|\frac{\sigma_0^2}{h_0}\right\|_K < \infty \quad \text{and} \quad \mathbb{E}\|\log h_0\|_K < \infty.$$

These conditions are similar to those of [19].

THEOREM 4.1. *Under the conditions* C.1–C.4, *the QMLE* $\hat{\theta}_n$ *is strongly consistent, that is,*

$$\hat{\boldsymbol{\theta}}_n \xrightarrow{\text{a.s.}} \boldsymbol{\theta}_0, \qquad n \to \infty.$$

PROOF.  In the first part we give the proof under the additional condition C.5 which allows one to apply Theorem 2.7 to $L_n/n$. In the second part we will indicate that C.5 is not needed. We have chosen to give two different proofs because the use of the uniform strong law of large numbers is intuitively more appealing than the proof without C.5.

*Part* 1. Assume C.1–C.5. First we show that $\hat{L}_n/n \xrightarrow{\text{a.s.}} L$ in $\mathbb{C}(K)$ as $n \to \infty$, where

$$L(\boldsymbol{\theta}) = -\frac{1}{2}\mathbb{E}\left(\frac{\sigma_0^2}{h_0(\boldsymbol{\theta})} + \log h_0(\boldsymbol{\theta})\right), \qquad \boldsymbol{\theta} \in K.$$

Second, we need to prove that $L$ is uniquely maximized at $\boldsymbol{\theta} = \boldsymbol{\theta}_0$. In the third step we show that the a.s. uniform convergence of $\hat{L}_n/n$ toward $L$ together with the fact that the limit $L$ has a unique maximum implies strong consistency.

(i) We first establish $L_n/n \xrightarrow{\text{a.s.}} L$ in $\mathbb{C}(K)$ by an application of Theorem 2.7. Proposition 3.12 shows that the sequence $(\ell_t) = -2^{-1}(X_t^2/h_t + \log h_t)$ of random elements with values in $\mathbb{C}(K)$ is of the form $(\ell_t) = (f(X_t, X_{t-1}, \dots))$, where $f$ is measurable, and hence, stationary ergodic (Proposition 2.5). Since $X_0 = \sigma_0 Z_0$ with $Z_0$ independent of $\sigma_0$ and $h_0$ and $\mathbb{E}Z_0^2 = 1$, assumption C.5 implies $\mathbb{E}\|X_0^2/h_0\|_K = \mathbb{E}\|\sigma_0^2/h_0\|_K < \infty$. Altogether, $\mathbb{E}\|\ell_0\|_K < \infty$, so that $L_n/n \xrightarrow{\text{a.s.}} L$ by Theorem 2.7. The property $\hat{L}_n/n \xrightarrow{\text{a.s.}} L$ follows if we can demonstrate $\|\hat{L}_n - L_n\|_K/n \xrightarrow{\text{a.s.}} 0$. Since $h_t, \hat{h}_t \geq \underline{g} > 0$, an application of the mean value theorem leads to $\|(\hat{h}_t)^{-1} - (h_t)^{-1}\|_K \leq \underline{g}^{-2}\|\hat{h}_t - h_t\|_K$ and $\|\log \hat{h}_t - \log h_t\|_K \leq \underline{g}^{-1}\|\hat{h}_t - h_t\|_K$. Thus, there exists $c > 0$ with

$$\|\hat{L}_n - L_n\|_K \leq c\sum_{t=1}^{n}(1 + X_t^2)\|\hat{h}_t - h_t\|_K \leq c\sum_{t=1}^{\infty}(1 + X_t^2)\|\hat{h}_t - h_t\|_K.$$

By Proposition 3.12, $\|\hat{h}_t - h_t\|_K \xrightarrow{\text{e.a.s.}} 0$ and by condition C.1 together with $\mathbb{E}Z_0^2 = 1$ and Lemma 2.2, we have $\mathbb{E}[\log^+(1 + X_0^2)] < \infty$. An application of Lemma 2.1 demonstrates $\sum_{t=1}^{\infty}(1 + X_t^2)\|\hat{h}_t - h_t\|_K < \infty$ a.s. Hence, $\|\hat{L}_n - L_n\|_K/n \xrightarrow{\text{a.s.}} 0$ and $\hat{L}_n/n \xrightarrow{\text{a.s.}} L$, as claimed.

(ii) For the uniqueness of the maximum of $L$ on $K$, we need to prove that $L(\boldsymbol{\theta}) < L(\boldsymbol{\theta}_0)$ for all $\boldsymbol{\theta} \in K \setminus \{\boldsymbol{\theta}_0\}$. Since $\mathbb{E}[\log \sigma_0^2]$ is finite and does not depend on the parameter $\boldsymbol{\theta}$, we can equivalently demonstrate that the function

$$Q(\boldsymbol{\theta}) = \mathbb{E}\left(\log \frac{\sigma_0^2}{h_0(\boldsymbol{\theta})} - \frac{X_0^2}{h_0(\boldsymbol{\theta})}\right) = \mathbb{E}\left(\log \frac{\sigma_0^2}{h_0(\boldsymbol{\theta})} - \frac{\sigma_0^2}{h_0(\boldsymbol{\theta})}\right), \qquad \boldsymbol{\theta} \in K,$$

is uniquely maximized at $\boldsymbol{\theta} = \boldsymbol{\theta}_0$. One can verify that $\log(x) - x \leq -1$ for all $x > 0$ with equality if and only if $x = 1$. Hence, $Q(\boldsymbol{\theta}) \leq -1 = Q(\boldsymbol{\theta}_0)$ with equality if and only if $\sigma_0^2/h_0(\boldsymbol{\theta}) \equiv 1$ a.s. By C.4, $\sigma_0^2/h_0(\boldsymbol{\theta}) \equiv 1$ if and only if $\boldsymbol{\theta} = \boldsymbol{\theta}_0$, which shows that $Q$ and $L$ are uniquely maximized at $\boldsymbol{\theta} = \boldsymbol{\theta}_0$.

(iii) Showing that (i) and (ii) imply strong consistency is accomplished by using standard arguments, which go back to [32].

*Part* 2. Without C.5 there is no longer uniform convergence of $L_n/n$ toward $L$. The proof of strong consistency rests on an argument by Pfanzagl [28]. By virtue of Proposition 3.12 and C.3, the function

$$\boldsymbol{\theta} \mapsto \ell_t(\boldsymbol{\theta}) = -\frac{1}{2}\left(\frac{X_t^2}{h_t(\boldsymbol{\theta})} + \log h_t(\boldsymbol{\theta})\right)$$

is continuous on $K$ with probability 1. Since, for every *fixed* $\boldsymbol{\theta} \in K$, the sequence $(\ell_t(\boldsymbol{\theta}))$ is stationary ergodic, one has that $n^{-1}\sum_{t=1}^{n}\ell_t(\boldsymbol{\theta}) \xrightarrow{\text{a.s.}} L(\boldsymbol{\theta}) = \mathbb{E}[\ell_0(\boldsymbol{\theta})]$ as $n \to \infty$ by an application of the ergodic theorem for real random variables; note that if $\mathbb{E}X_0^2 = \infty$, the latter limit can take the value $-\infty$ at certain points $\boldsymbol{\theta}$, but $h_0(\boldsymbol{\theta}) \geq \underline{g} > 0$ guarantees $L(\boldsymbol{\theta}) < +\infty$ for all $\boldsymbol{\theta} \in K$. Therefore, we can use exactly the same arguments as given in the proof of Lemma 3.11 of [28] in order to show that the function $L$ is upper semicontinuous on $K$ and $\limsup_{n \to \infty} \sup_{\boldsymbol{\theta} \in K'} L_n(\boldsymbol{\theta})/n \leq \sup_{\boldsymbol{\theta} \in K'} L(\boldsymbol{\theta})$ with probability 1 for any compact subset $K' \subseteq K$. Since C.1–C.4 imply $\|\hat{L}_n - L_n\|_K/n \xrightarrow{\text{a.s.}} 0$, the inequality $\limsup_{n \to \infty} \sup_{\boldsymbol{\theta} \in K'} \hat{L}_n(\boldsymbol{\theta})/n \leq \sup_{\boldsymbol{\theta} \in K'} L(\boldsymbol{\theta})$ a.s. is valid also. Because an upper semicontinuous function attains its maximum on compact sets, one can demonstrate $\hat{\boldsymbol{\theta}}_n \xrightarrow{\text{a.s.}} \boldsymbol{\theta}_0$, similarly to step (iii) of the proof of Part 1. □

## 5. Examples: consistency.
For the purpose of illustration, we apply Theorem 4.1 to EGARCH and AGARCH$(p, q)$. Our task will be to define the set $K$ in an appropriate way and to verify the conditions C.1–C.4.

5.1. *EGARCH.* Subsume the EGARCH parameters $\alpha$, $\beta$, $\gamma$ and $\delta$ into $\boldsymbol{\theta} = (\alpha, \beta, \gamma, \delta)^T$ and denote by $\boldsymbol{\theta}_0 = (\alpha_0, \beta_0, \gamma_0, \delta_0)^T$ the true parameter vector. As discussed in Example 3.8, we suppose $0 \le \beta_0 < 1$, $\delta_0 \ge |\gamma_0|$. The parameter space is of the form $\Theta = \mathbb{R} \times [0, 1) \times D_E$, where

$$D_E = \{(\gamma, \delta)^T \in \mathbb{R}^2 | \gamma \in \mathbb{R}, \delta \ge |\gamma|\}.$$

The compact set $K \subset \Theta$ will be defined below. The restriction $0 \le \beta_0 < 1$ guarantees that the SRE (3.7) has a unique stationary ergodic solution $(\log \sigma_t^2)$. From the almost sure representation (3.8) of this solution, one recognizes $\mathbb{E}(\log \sigma_0^2) < \infty$, which establishes C.1.

Rather than checking condition C.2 of Theorem 4.1, we *directly* verify its consequences which were used in the proof of the latter theorem, namely, the fact that a stationary ergodic sequence $(h_t)$ of random elements with values in $\mathbb{C}(K, [0, \infty))$ can be defined such that

(5.1)    $h_t$ is $\mathcal{F}_{t-1}$-measurable   and   $h_t(\boldsymbol{\theta}_0) = \sigma_t^2$   a.s. for every $t$,

(5.2)
$$X_t^2 \|\hat{h}_t^{-1} - h_t^{-1}\|_K \xrightarrow{\text{e.a.s.}} 0 \quad \text{and} \quad \|\log \hat{h}_t - \log h_t\|_K \xrightarrow{\text{e.a.s.}} 0$$

$$\text{as } t \to \infty.$$

To this end, we consider the SRE

(5.3)                    $$\log s_{t+1} = \Phi_t(\log s_t), \qquad t \in \mathbb{Z},$$

where

$$[\Phi_t(s)](\boldsymbol{\theta}) = \alpha + \beta \log s(\boldsymbol{\theta}) + (\gamma X_t + \delta |X_t|) \exp(-s(\boldsymbol{\theta})/2), \qquad \boldsymbol{\theta} \in K.$$

Observing that $(\gamma X_t + \delta |X_t|) \exp(-s(\boldsymbol{\theta})/2) \ge 0$, we find that for any constant function $s = \log \varsigma_0^2$ and $\varepsilon > 0$ the forward (and backward) iterates associated with (5.3) satisfy

$$\log \hat{h}_t(\boldsymbol{\theta}) = [\Phi_{t-1}^{(t)}(\log \varsigma_0^2)](\boldsymbol{\theta}) \ge \alpha(1 + \beta + \cdots + \beta^{t-1}) + \beta^t \log \varsigma_0^2$$

$$\ge \alpha(1 - \beta)^{-1} - \varepsilon$$

for large enough $t$. Therefore, we may suppose without loss of generality that the SRE (5.3) lives on the subset $\mathbb{C}(K, [m - \varepsilon, \infty))$ with $m = \inf_{\boldsymbol{\theta} \in K} \alpha(1 - \beta)^{-1}$. By comparison with the derivation of (3.22), one recognizes that

$$\Lambda(\Phi_0) = \sup_{\boldsymbol{\theta} \in K} \lambda_\varepsilon(\boldsymbol{\theta}) = \|\lambda_\varepsilon\|_K,$$

where

(5.4)   $\lambda_\varepsilon(\boldsymbol{\theta}) = \max(\beta, 2^{-1} \exp(-(m - \varepsilon)/2)(\gamma X_0 + \delta |X_0|) - \beta), \qquad \boldsymbol{\theta} \in K.$

Since $\mathbb{E}(\log \|\lambda_\varepsilon\|_K) \to \mathbb{E}(\log \|\lambda_0\|_K)$ as $\varepsilon \downarrow 0$, the condition $\mathbb{E}(\log \|\lambda_0\|_K) < 0$ implies $\mathbb{E}(\log \|\lambda_\varepsilon\|_K) < 0$ for small enough $\varepsilon$. Hence, if we assume $K$ is

chosen such that $\mathbb{E}(\log\|\lambda_0\|_K) < 0$, then the sequence $(\Lambda(\Phi_t))$ obeys the conditions of Theorem 2.8, and thus, the SRE (5.3) (on $\mathbb{C}(K, [m - \varepsilon, \infty))$, where $\varepsilon > 0$ small) admits a unique stationary solution $(\log h_t)$, which is ergodic. Moreover, $\log h_t$ is $\mathcal{F}_{t-1}$-measurable and the property $h_t(\boldsymbol{\theta}) = h_t(\boldsymbol{\theta}_0)$ a.s. follows from similar arguments as in the proof of Proposition 3.12. Thus, (5.1) is established. As regards the verification of (5.2), note that $\|\log \hat{h}_t - \log h_t\|_K \overset{\text{e.a.s.}}{\longrightarrow} 0$ by Theorem 2.8. The mean value theorem applied to the function $e^{-x}$ together with the facts that $\log \hat{h}_t \geq m - \varepsilon$ for $t$ large and $\log h_t \geq m$ for all $t$ yields a constant $c > 0$ with

$$X_t^2 \|(\hat{h}_t)^{-1} - (h_t)^{-1}\|_K = X_t^2 \|\exp(-\log \hat{h}_t) - \exp(-\log h_t)\|_K$$
$$\leq c X_t^2 \|\log \hat{h}_t - \log h_t\| \overset{\text{e.a.s.}}{\longrightarrow} 0.$$

The latter limit relation is a consequence of $\mathbb{E}(\log^+ X_0^2) < \infty$, implied by $\mathbb{E}(\log \sigma_0^2) < \infty$ and $\mathbb{E}(Z_0^2) = 1$, and an application of Lemma 2.1. This completes the demonstration of the limit relations (5.2).

Since condition C.3 is automatically satisfied for the EGARCH process, we are left with the verification of the identifiability condition C.4. Before we start, we impose the mild technical assumptions that $(\gamma_0, \delta_0) \neq (0, 0)$ and that the distribution of $Z_0$ is not concentrated at two points. Note that the case $\gamma_0 = \delta_0 = 0$ would lead to an identifiability problem because then $\log \sigma_t^2 = \alpha_0(1 - \beta_0)^{-1}$ by representation (3.8), which implies that there are infinitely many parameters leading to the identical model. Observe that $h_0(\boldsymbol{\theta}) = \sigma_0^2$ a.s. is equivalent to $\log h_t(\boldsymbol{\theta}) = \log h_t(\boldsymbol{\theta}_0)$ a.s. for all $t \in \mathbb{Z}$ because of the stationarity of $(\log h_t - \sigma_t^2)$ and the property $\log h_0(\boldsymbol{\theta}_0) = \log \sigma_0^2$. We now show the nontrivial implication $\log h_t(\boldsymbol{\theta}) = \log h_t(\boldsymbol{\theta}_0) \Rightarrow \boldsymbol{\theta} = \boldsymbol{\theta}_0$. Replacing $\log h_t$ by $\Phi_{t-1}(\log h_{t-1})$ in the identity $\log h_t(\boldsymbol{\theta}) = \log h_t(\boldsymbol{\theta}_0)$ and accounting for $h_{t-1}(\boldsymbol{\theta}) = \sigma_{t-1}^2$, we obtain

$$(\alpha - \alpha_0) + (\beta - \beta_0)\log \sigma_{t-1}^2 + \{(\gamma - \gamma_0)Z_{t-1} + (\delta - \delta_0)|Z_{t-1}|\} = 0 \qquad \text{a.s.}$$

If $\beta \neq \beta_0$, the random variable $\log \sigma_{t-1}^2$ would at the same time be a measurable function of $Z_{t-1}$ and independent of $Z_{t-1}$. This implies that $\log \sigma_{t-1}^2$ is deterministic. However, taking the variance in (3.8) gives $\text{Var}(\log \sigma_{t-1}^2) = \sum_{k=0}^{\infty} \beta_0^{2k} \text{Var}(\gamma_0 Z_0 + \delta_0 |Z_0|) > 0$ since the facts that $(\gamma_0, \delta_0) \neq (0, 0)$ and that the distribution of $Z_0$ is not concentrated at two points imply $\text{Var}(\gamma_0 Z_0 + \delta_0 |Z_0|) > 0$. To avoid this contradiction, necessarily $\beta = \beta_0$. Furthermore, if $(\alpha - \alpha_0, \gamma - \gamma_0, \delta - \delta_0) \neq (0, 0, 0)$, there are three distinct cases concerning the zeros of the function $f(z) = (\alpha - \alpha_0) + (\gamma - \gamma_0)z + (\delta - \delta_0)|z|$:

(1) $f$ has less than or equal to two zeros.
(2) $f \equiv 0$ on $[0, \infty)$ and $f \neq 0$ on $(-\infty, 0)$.
(3) $f \neq 0$ on $(0, \infty)$ and $f \equiv 0$ on $(-\infty, 0]$.

Using this observation and the facts that $\mathbb{E}Z_0 = 0$, $\mathbb{E}Z_0^2 = 1$ and that the distribution of $Z_0$ is not concentrated at two points, we conclude that $f(Z_{t-1}) = 0$ a.s. if and only if $(\alpha, \beta, \gamma) = (\alpha_0, \gamma_0, \delta_0)$. Thus, $\boldsymbol{\theta} = \boldsymbol{\theta}_0$, which concludes the verification of C.4. Therefore, the strong consistency of the QMLE in EGARCH can be established by an application of Theorem 4.1.

THEOREM 5.1. *Let* $((X_t, \sigma_t))$ *be a stationary EGARCH process with parameters* $\boldsymbol{\theta}_0 = (\alpha_0, \beta_0, \gamma_0, \delta_0)^T$ *such that* $0 \leq \beta_0 < 1$ *and* $\delta_0 \geq |\gamma_0|$ *and* $(\gamma_0, \delta_0) \neq (0, 0)$. *Suppose the distribution of* $Z_0$ *is not concentrated at two points. Let* $K \subset \mathbb{R} \times [0, 1) \times D_E$ *be a compact set with* $\boldsymbol{\theta}_0 \in K$ *and such that*

$$\mathbb{E}(\log \|\lambda_0\|_K) < 0,$$

*where* $\lambda_0$ *is given by* (5.4) *with* $m = \inf_{\boldsymbol{\theta} \in K} \alpha(1-\beta)^{-1}$. *Then the QMLE* (4.2) *is strongly consistent.*

5.2. AGARCH$(p, q)$. Set $\boldsymbol{\theta} = (\alpha_0, \alpha_1, \ldots, \alpha_p, \beta_1, \ldots, \beta_q, \gamma)^T \in \mathbb{R}^{p+q+1}$ and denote the true parameter vector of AGARCH$(p, q)$ by $\boldsymbol{\theta}_0 = (\alpha_0^\circ, \alpha_1^\circ, \ldots, \alpha_p^\circ, \beta_1^\circ, \ldots, \beta_q^\circ, \gamma^\circ)^T$. We suppose that $\boldsymbol{\theta}_0$ admits a unique stationary ergodic solution $((X_t, \sigma_t))$ to the AGARCH$(p, q)$ equations, which is equivalent to a strictly negative top Lyapunov exponent of the associated matrix sequence $(\mathbf{A}_t)$ given by (3.12) with $\alpha_i = \alpha_i^\circ$, $\beta_j = \beta_j^\circ$ and $\gamma = \gamma^\circ$; see Theorem 3.5. Moreover, suppose $\alpha_i^\circ > 0$ for some $i = 1, \ldots, p$ because otherwise the constant sequence $\sigma_t^2 = \alpha_0^\circ(1 - \sum_{j=1}^q \beta_j^\circ)^{-1}$ is the unique stationary solution of (3.9), which would imply that one cannot discriminate between $\alpha_0^\circ$ and the $\beta_j^\circ$'s (nonidentifiability). Another necessary restriction is $(\alpha_p^\circ, \beta_q^\circ) \neq (0, 0)$. Let $K$ be a compact subset of $(0, \infty) \times [0, \infty)^p \times B \times [-1, 1]$ containing the true parameter $\boldsymbol{\theta}_0$, where

$$(5.5) \qquad B = \left\{ (\beta_1, \ldots, \beta_q)^T \in [0, 1)^q \,\Big|\, \sum_{j=1}^q \beta_j < 1 \right\}.$$

We now verify C.1–C.4.

In Remark 3.6 we have shown the existence of a $\tilde{q} > 0$ with $\mathbb{E}\sigma_0^{2\tilde{q}} < \infty$. Thus, $\mathbb{E}(\log^+ \sigma_0^2) < \infty$ is valid.

As regards C.2, note that

$$[\Phi_t(\mathbf{s})](\boldsymbol{\theta}) = \left( \alpha_0 + \sum_{i=1}^p \alpha_i(|X_{t+1-i}| - \gamma X_{t+1-i})^2 \right) \mathbf{e}_1 + \mathbf{C}(\boldsymbol{\theta})\mathbf{s},$$

for $\mathbf{s} \in \mathbb{C}(K, [0, \infty)^q)$, where

$$\mathbf{e}_1 = (1, 0, \ldots, 0)^T \in \mathbb{R}^q,$$

$$\mathbf{C}(\boldsymbol{\theta}) = \begin{pmatrix} \beta_1 & \beta_2 & \cdots & \cdots & \beta_{q-1} & \beta_q \\ 1 & 0 & \cdots & \cdots & 0 & 0 \\ 0 & 1 & 0 & \cdots & \cdots & 0 \\ \vdots & \ddots & \ddots & \ddots & \vdots & \vdots \\ 0 & \cdots & 0 & 1 & 0 & 0 \\ 0 & \cdots & \cdots & 0 & 1 & 0 \end{pmatrix} \in \mathbb{R}^{q \times q};$$

see Example 3.11. Analogously to the demonstration of invertibility in AGARCH, $\Lambda(\Phi_0^{(r)}) \leq \Lambda(\mathbf{C}^r)$, $r \in \mathbb{N}$. Note that $\bar{\beta} := \sup_{\boldsymbol{\theta} \in K}(\sum_{j=1}^q \beta_j) < 1$ since $K$ is compact. Pointwise application of the inequality (3.25) for each $\boldsymbol{\theta} \in K$ yields for any $\mathbf{s} \in \mathbb{C}(K, [0,\infty)^q)$ the bound

$$|(\mathbf{C}^r \mathbf{s})(\boldsymbol{\theta})| \leq \mathrm{const}\, \bar{\beta}^{r/q} |\mathbf{s}(\boldsymbol{\theta})|,$$

where the constant does not depend on $\mathbf{s}$ and $\boldsymbol{\theta}$. Taking the supremum on both sides of the latter bound, one obtains

$$\|\mathbf{C}^r \mathbf{s}\|_K \leq \mathrm{const}\, \bar{\beta}^{r/q} \|\mathbf{s}\|_K, \qquad r \geq 0,$$

showing that $\Lambda(\mathbf{C}^r) \leq \mathrm{const}\, \bar{\beta}^{r/q} \to 0$ as $r \to \infty$. Therefore, conditions 1 and 2 of Proposition 3.12 are verified. Hence, C.2 holds and $\mathbf{h}_t = (h_t, \ldots, h_{t-q+1})^T$ is properly defined.

C.3 being obviously satisfied, we turn to the identifiability condition C.4. We split our arguments into a series of lemmas. First we derive an almost sure representation of $h_t$, similarly to [2].

Lemma 5.2. *The following almost sure representation for $h_t$ is valid:*

$$(5.6) \qquad h_t(\boldsymbol{\theta}) = \xi_0(\boldsymbol{\theta}) + \sum_{\ell=1}^\infty \xi_\ell(\boldsymbol{\theta})(|X_{t-\ell}| - \gamma X_{t-\ell})^2, \qquad \boldsymbol{\theta} \in K,$$

*with the sequence $(\xi_\ell(\boldsymbol{\theta}))_{\ell \in \mathbb{N}}$ given by*

$$(5.7) \qquad \begin{aligned} \xi_0(\boldsymbol{\theta}) &= \frac{\alpha_0}{b_{\boldsymbol{\theta}}(1)} = \alpha_0\left(1 - \sum_{j=1}^q \beta_j\right)^{-1} \quad and \\ \sum_{\ell=1}^\infty \xi_\ell(\boldsymbol{\theta})z^\ell &= \frac{a_{\boldsymbol{\theta}}(z)}{b_{\boldsymbol{\theta}}(z)}, \qquad |z| \leq 1, \end{aligned}$$

*where $a_{\boldsymbol{\theta}}(z) = \sum_{i=1}^p \alpha_i z^i$ and $b_{\boldsymbol{\theta}}(z) = 1 - \sum_{j=1}^q \beta_j z^j$.*

Proof. The proof rests on the observation that $(h_t(\boldsymbol{\theta}))$ obeys an ARMA$(p,q)$ equation, that is,

$$h_t(\boldsymbol{\theta}) = \alpha_0 + \sum_{i=1}^p \alpha_i(|X_{t-i}| - \gamma X_{t-i})^2 + \sum_{j=1}^q \beta_j h_{t-j}(\boldsymbol{\theta}), \qquad t \in \mathbb{Z},$$

or shorter, in backshift operator notation,

$$(5.8) \qquad b_{\boldsymbol{\theta}}(B)h_t(\boldsymbol{\theta}) = \alpha_0 + a_{\boldsymbol{\theta}}(B)(|X_t| - \gamma X_t)^2, \qquad t \in \mathbb{Z}.$$

Since $\bar{\beta} < 1$, the zeros of the polynomial $b_{\boldsymbol{\theta}}(z)$ lie outside the unit disc; indeed, if $|z| < \bar{\beta}^{-1/q}$, then $|b_{\boldsymbol{\theta}}(z)| \geq 1 - \sum_{j=1}^{q} \beta_j |z|^j > 1 - \bar{\beta}(\bar{\beta}^{-1/q})^q = 0$. This suggests the a.s. representation (5.6), as can be seen from a comparison with Section 3.2 of [11]. To prove this a.s. representation, we cannot, however, directly apply Proposition 3.1.2 in [11] because the "innovations" $(|X_t| - \gamma X_t)^2$ may have an infinite second moment. One possible way to validate the a.s. representation is to show that the right-hand side of (5.6) is well-defined, continuous on $K$ a.s., and that it obeys (5.8) regarded as a difference equation on $\mathbb{C}(K)$. Then, since the SRE for $(\mathbf{h}_t)$ admits a *unique* stationary ergodic solution $(\mathbf{h}_t)$, the right-hand side of (5.6) must coincide with $h_t$. To show the three assertions mentioned before, first note that there are $0 < \eta < 1$ and $c > 0$ with $|\xi_\ell(\boldsymbol{\theta})| \leq c\eta^\ell$ for all $\ell \geq 1$ and $\boldsymbol{\theta} \in K$ [apply the Cauchy inequalities to the complex function $1/b_{\boldsymbol{\theta}}(z)$] and, thus, $\xi_\ell \xrightarrow{\text{e.a.s.}} 0$ in $\mathbb{C}(K)$ as $\ell \to \infty$. Since $\xi_\ell$ is continuous on $K$ and $\mathbb{E}[(|X_0| - \gamma^\circ X_0)^{2\tilde{q}}] < \infty$ for a $\tilde{q} > 0$ (see Remark 3.6), the series $\sum_{\ell=1}^{\infty} \xi_\ell(\boldsymbol{\theta})(|X_{t-\ell}| - \gamma X_{t-\ell})^2$ converges absolutely a.s. in $\mathbb{C}(K)$ by virtue of Lemma 2.1. Hence, (5.6) is continuous a.s. Eventually it is an elementary exercise to prove that the a.s. representation for $(h_t)$ obeys (5.8). This completes the proof of the lemma. $\quad \square$

The next lemma is concerned with the identifiability of the parameter $\gamma$.

LEMMA 5.3. *Suppose that the distribution of $Z_0$ is not concentrated at two points. Then for any $\boldsymbol{\theta} \in K$, the relation $h_0(\boldsymbol{\theta}) = h_0(\boldsymbol{\theta}_0)$ a.s. implies $\gamma = \gamma^\circ$.*

PROOF. Note that $h_0(\boldsymbol{\theta}) = h_0(\boldsymbol{\theta}_0)$ a.s. is equivalent to $h_k(\boldsymbol{\theta}) = h_k(\boldsymbol{\theta}_0)$ for any $k$, in particular, $k = \max(p, q)$. We rewrite $h_k(\boldsymbol{\theta}) = h_k(\boldsymbol{\theta}_0)$ a.s. as

$$(5.9) \qquad (\alpha_0^\circ - \alpha_0) + \sum_{i=1}^{k} Y_{k-i}\sigma_{k-i}^2 = 0 \qquad \text{a.s.,}$$

where $Y_{k-i} = \alpha_i^\circ(|Z_{k-i}| - \gamma^\circ Z_{k-i})^2 - \alpha_i(|Z_{k-i}| - \gamma Z_{k-i})^2 + (\beta_i^\circ - \beta_i)$. Introduce $k^* = \min(i \in [1, p] \mid \alpha_i^\circ > 0)$. Then by repeatedly expressing each term $\sigma_{t-j}^2$, $j = 1, \ldots, (k^* - 1)$, in the equation

$$\sigma_t^2 = \alpha_0^\circ + \sum_{i=k^*}^{p} \alpha_i^\circ(|X_{t-i}| - \gamma^\circ X_{t-i})^2 + \sum_{j=1}^{q} \beta_j^\circ \sigma_{t-j}^2$$

by past observations and past squared volatilities, one sees that $\sigma_t^2$ can be written as a function of $\{(|X_{t-k^*}| - \gamma^\circ X_{t-k^*})^2, (|X_{t-k^*-1}| - \gamma^\circ X_{t-k^*-1})^2, \ldots;$

$\sigma_{t-k^*}^2, \sigma_{t-k^*-1}^2, \ldots\}$, and consequently, $\sigma_t^2$ is $\mathcal{F}_{t-k^*}$-measurable. Relation (5.9) together with $\sigma_{k-1}^2 \geq \alpha_0^\circ > 0$ implies that $Y_{k-1}$ is a function of $\sigma_0^2, \ldots, \sigma_{k-1}^2$ and consequently, $\mathcal{F}_{k-k^*-1}$-measurable. Since $Y_{k-1}$ is at the same time independent of $\mathcal{F}_{k-k^*-1}$, it must be degenerate. With the identical arguments, $Y_{k-2}, \ldots, Y_{k-k^*}$ are degenerate. The degeneracy of $Y_{k-k^*}$ means that

$$\alpha_{k^*}^\circ (|Z_{k-k^*}| - \gamma^\circ Z_{k-k^*})^2 - \alpha_{k^*}(|Z_{k-k^*}| - \gamma Z_{k-k^*})^2 = \tilde{c}$$

for a certain constant $\tilde{c}$. Note that a.s. on the sets $\{Z_{k-k^*} \geq 0\}$ and $\{Z_{k-k^*} < 0\}$,

$$(\alpha_{k^*}^\circ (1-\gamma^\circ)^2 - \alpha_{k^*}(1-\gamma)^2) Z_{k-k^*}^2 = \tilde{c}$$

and

$$(\alpha_{k^*}^\circ (1+\gamma^\circ)^2 - \alpha_{k^*}(1+\gamma)^2) Z_{k-k^*}^2 = \tilde{c},$$

respectively. Because the distribution of $Z_{k-k^*}$ is not concentrated at two points, these two equations can only be jointly satisfied if $\tilde{c} = 0$. Since $\mathbb{E}Z_{k-k^*} = 0$ and $\mathbb{E}Z_{k-k^*}^2 = 1$, the distribution of $Z_{k-k^*}$ has positive mass both on the negative and the positive real line, which implies $\alpha_{k^*}^\circ (1-\gamma^\circ)^2 = \alpha_{k^*}(1-\gamma)^2$ and $\alpha_{k^*}^\circ (1+\gamma^\circ)^2 = \alpha_{k^*}(1+\gamma)^2$. From these two equations together with $\alpha_{k^*}^\circ > 0$, we conclude $\alpha_{k^*} > 0$. Subtracting and adding the latter two equations yields $\gamma = \gamma^\circ \alpha_{k^*}^\circ / \alpha_{k^*}$ and $1 + \gamma^2 = (1 + (\gamma^\circ)^2) \alpha_{k^*}^\circ / \alpha_{k^*}$. Because of the constraint $|\gamma| \leq 1$, we can conclude $\gamma = \gamma^\circ$ (and $\alpha_{k^*} = \alpha_{k^*}^\circ$), which completes the proof. □

Eventually we establish the identifiability condition C.4.

LEMMA 5.4. *Suppose that the distribution of $Z_0$ is not concentrated at two points and that the polynomials $a_{\boldsymbol{\theta}_0}$ and $b_{\boldsymbol{\theta}_0}$ defined in Lemma 5.2 do not have any common zeros. Then for any $\boldsymbol{\theta} \in K$,*

$$h_0(\boldsymbol{\theta}) = \sigma_0^2 \quad \text{if and only if} \quad \boldsymbol{\theta} = \boldsymbol{\theta}_0.$$

PROOF. We have shown in Lemma 5.3 that $h_0(\boldsymbol{\theta}) = \sigma_0^2$ implies $\gamma = \gamma^\circ$, so that in consideration of (5.6), the relation $h_0(\boldsymbol{\theta}) = h_0(\boldsymbol{\theta}_0)$ becomes

$$\xi_0(\boldsymbol{\theta}) - \xi_0(\boldsymbol{\theta}_0) + \sum_{\ell=1}^\infty (\xi_\ell(\boldsymbol{\theta}) - \xi_\ell(\boldsymbol{\theta}_0))(|X_{-\ell}| - \gamma^\circ X_{-\ell})^2 \equiv 0.$$

We first show $\xi_\ell(\boldsymbol{\theta}) - \xi_\ell(\boldsymbol{\theta}_0) = 0$ for all $\ell \in \mathbb{N}$ by contradiction. Denote by $\ell^* \geq 1$ the smallest integer $\ell \geq 1$ with $\delta_\ell := \xi_\ell(\boldsymbol{\theta}) - \xi_\ell(\boldsymbol{\theta}_0) \neq 0$. Since $\sigma_{-\ell^*}^2 \geq \alpha_0^\circ > 0$, we have that

$(|Z_{-\ell^*}| - \gamma^\circ Z_{-\ell^*})^2$

$$= \left( \xi_0(\boldsymbol{\theta}_0) - \xi_0(\boldsymbol{\theta}) + \sum_{\ell=\ell^*+1}^\infty (\xi_\ell(\boldsymbol{\theta}_0) - \xi_\ell(\boldsymbol{\theta}))(|X_{-\ell}| - \gamma^\circ X_{-\ell})^2 \right) \Big/ (\delta_{\ell^*} \sigma_{-\ell^*}^2)$$

is at the same time $\mathcal{F}_{-\ell^*-1}$-measurable and independent of $\mathcal{F}_{-\ell^*-1}$, which is only possible if $(|Z_{-\ell^*}| - \gamma^\circ Z_{-\ell^*})^2$ is degenerate. However, from the assumption that the distribution of $Z_0$ is not concentrated at two points, it follows that $(|Z_{-\ell^*}| - \gamma^\circ Z_{-\ell^*})^2$ cannot be degenerate, that is, the desired contradiction. Using $\xi_\ell(\boldsymbol{\theta}) = \xi_\ell(\boldsymbol{\theta}_0)$ in (5.7), we conclude $a_{\boldsymbol{\theta}}(z)/b_{\boldsymbol{\theta}}(z) = a_{\boldsymbol{\theta}_0}(z)/b_{\boldsymbol{\theta}_0}(z)$. Write $a_{\boldsymbol{\theta}}(z) = r(z)a_{\boldsymbol{\theta}_0}(z)$ and $b_{\boldsymbol{\theta}}(z) = r(z)b_{\boldsymbol{\theta}_0}(z)$. The rational function $r(z)$ does not have any pole because otherwise $a_{\boldsymbol{\theta}_0}(z)$ and $b_{\boldsymbol{\theta}_0}(z)$ would have a common zero. Hence, $r(z)$ is a polynomial. The degree of $r$ is zero, because otherwise either $a_{\boldsymbol{\theta}}(z)$ or $b_{\boldsymbol{\theta}}(z)$ would have degree strictly greater than $p$ or $q$, respectively, since $(\alpha_p^\circ, \beta_q^\circ) \neq (0, 0)$. Finally, $r \equiv 1$ because the constants in the polynomials $b_{\boldsymbol{\theta}}$ and $b_{\boldsymbol{\theta}_0}$ are 1. Hence, $a_{\boldsymbol{\theta}} = a_{\boldsymbol{\theta}_0}$ and $b_{\boldsymbol{\theta}} = b_{\boldsymbol{\theta}_0}$, which gives $\boldsymbol{\theta} = \boldsymbol{\theta}_0$ and concludes the proof. $\square$

Now an application of Theorem 4.1 yields strong consistency of the QMLE in AGARCH$(p,q)$. This result generalizes Theorem 4.1 of [2].

THEOREM 5.5. *Let $(X_t)$ be a stationary* AGARCH$(p,q)$ *process with true parameters $\boldsymbol{\theta}_0 = (\alpha_0^\circ, \alpha_1^\circ, \ldots, \alpha_p^\circ, \beta_1^\circ, \ldots, \beta_q^\circ, \gamma^\circ)^T$ such that we have the following:*

1. $\alpha_i^\circ > 0$ *for some $i > 0$ and $(\alpha_p^\circ, \beta_q^\circ) \neq (0, 0)$.*
2. *The polynomials $a^\circ(z) = \sum_{i=1}^p \alpha_i^\circ z^i$ and $b^\circ(z) = 1 - \sum_{j=1}^q \beta_j^\circ z^j$ do not have any common zeros.*

*Suppose that the distribution of $Z_0$ is not concentrated at two points and let $K \subset (0, \infty) \times [0, \infty)^p \times B \times [-1, 1]$ be compact and contain $\boldsymbol{\theta}_0 \in K$, where $B$ is given in (5.5). Then the QMLE (4.2) is strongly consistent.*

**6. The first and second derivatives of the functions $h_t$ and $\hat{h}_t$.** For establishing the asymptotic normality of the QMLE, it is essential to understand the limit behavior of the sequences of functions $(\hat{h}'_t)_{t \in \mathbb{N}}$ and $(\hat{h}''_t)_{t \in \mathbb{N}}$ and to study the differentiability properties of $h_t$. The interpretation of the arising problems in terms of SRE's turns out to be fruitful once again. We first derive SRE's for the first and second derivatives of $\hat{\mathbf{h}}_t$ and in a second step we construct stationary approximations of $(\hat{\mathbf{h}}'_t)_{t \in \mathbb{N}}$ and $(\hat{\mathbf{h}}''_t)_{t \in \mathbb{N}}$, which turn out to coincide with $(\mathbf{h}'_t)_{t \in \mathbb{N}}$ and $(\mathbf{h}''_t)_{t \in \mathbb{N}}$, respectively.

For a notationally tractable representation of the SRE's to be derived, we introduce maps

$$\varphi_t: K \times [0, \infty)^q \to K \times [0, \infty)^q, \qquad (\boldsymbol{\theta}, \mathbf{u}) \mapsto (\boldsymbol{\theta}, (g_{\boldsymbol{\theta}}(\mathbf{X}_t, \mathbf{u}), u_1, \ldots, u_{q-1})),$$
$$t \in \mathbb{Z},$$

$$p_2: K \times [0, \infty)^q \to [0, \infty)^q, \qquad (\boldsymbol{\theta}, \mathbf{u}) \mapsto \mathbf{u},$$

and set $\psi_{t,r} = p_2 \circ \varphi_t^{(r)}$ for fixed $r \geq 1$. Observe that $[\Phi_t^{(r)}(\mathbf{s})](\boldsymbol{\theta}) = \psi_{t,r}(\boldsymbol{\theta}, \mathbf{s}(\boldsymbol{\theta}))$ for every $\mathbf{s} \in \mathbb{C}(K, [0, \infty)^q)$. In this (and only this) section, we work under the convention that the first- and second-order partial derivatives of a function $\mathbf{f} = (f_1, \ldots, f_m)^T : U \subset \mathbb{R}^n \to \mathbb{R}^m$ are written as vectors,

$$\mathbf{f}' = (\partial_1^1 \mathbf{f}, \partial_1^2 \mathbf{f}, \ldots, \partial_1^k \mathbf{f}, \ldots, \ldots, \partial_m^1 \mathbf{f}, \ldots, \partial_m^n \mathbf{f})^T,$$

$$\mathbf{f}'' = (\partial_1^{1,1} \mathbf{f}, \ldots, \partial_1^{1,n} \mathbf{f}, \ldots, \partial_1^{n,1} \mathbf{f}, \ldots, \partial_1^{n,n} \mathbf{f},$$
$$\partial_2^{1,1} \mathbf{f}, \ldots, \partial_2^{n,n} \mathbf{f}, \ldots, \partial_m^{1,1} \mathbf{f}, \ldots, \partial_m^{n,n} \mathbf{f})^T,$$

where $\partial_j^k \mathbf{f} := (\partial f_j)/(\partial x_k)$ and $\partial_j^{k_1,k_2} \mathbf{f} := (\partial^2 f_j)/(\partial x_{k_1} \partial x_{k_2})$. In what follows, we suppose that the following regularity conditions hold:

D.1 The conditions of Proposition 3.12 are satisfied with a compact set $K \subset \mathbb{R}^d$, which is contained in the interior of the parameter space $\Theta$. Suppose that $K$ coincides with the closure of its (open) interior. The function $(\boldsymbol{\theta}, s) \mapsto g_{\boldsymbol{\theta}}(\mathbf{x}, \mathbf{s})$ on $K \times [0, \infty)^q$ is continuously differentiable for every fixed $\mathbf{x} \in \mathbb{R}^p$.

D.2 For all $j \in \{1, \ldots, q\}$ and $k \in \{1, \ldots, d+q\}$

$$(6.1) \qquad \mathbb{E}\left[\log^+\left(\sup_{\boldsymbol{\theta} \in K} |\partial_j^k \psi_{0,1}(\boldsymbol{\theta}, \mathbf{h}_0(\boldsymbol{\theta}))|\right)\right] < \infty.$$

Moreover, there exist a stationary sequence $(\bar{C}_1(t))$ with $\mathbb{E}[\log^+ \bar{C}_1(0)] < \infty$ and $\kappa \in (0, 1]$ such that

$$(6.2) \qquad \sup_{\boldsymbol{\theta} \in K} |\partial_j^k \psi_{t,1}(\boldsymbol{\theta}, \mathbf{u}) - \partial_j^k \psi_{t,1}(\boldsymbol{\theta}, \tilde{\mathbf{u}})| \leq \bar{C}_1(t) |\mathbf{u} - \tilde{\mathbf{u}}|^\kappa, \qquad \mathbf{u}, \tilde{\mathbf{u}} \in [0, \infty)^q,$$

for every $j \in \{1, \ldots, q\}$, $k \in \{1, \ldots, d+q\}$ and $t \in \mathbb{Z}$.

We recall that the sequence $(\hat{\mathbf{h}}_t)_{t \in \mathbb{N}}$ of random elements with values in $\mathbb{C}(K, [0, \infty)^q)$ is a solution of the SRE $\mathbf{s}_{t+1} = \Phi_t(\mathbf{s}_t)$, $t \in \mathbb{N}$, on $\mathbb{C}(K, [0, \infty)^q)$; see Section 3.3. Differentiating with respect to $\boldsymbol{\theta}$ on both sides of $\hat{\mathbf{h}}_{t+1}(\boldsymbol{\theta}) = \Phi_t(\hat{\mathbf{h}}_t(\boldsymbol{\theta})) = \psi_{t,1}(\boldsymbol{\theta}, \hat{\mathbf{h}}_t(\boldsymbol{\theta}))$ yields

$$(6.3) \qquad \begin{aligned} \partial_j^k \hat{\mathbf{h}}_{t+1}(\boldsymbol{\theta}) &= \partial_j^k \psi_{t,1}(\boldsymbol{\theta}, \hat{\mathbf{h}}_t(\boldsymbol{\theta})) \\ &\quad + \sum_{i=1}^q \partial_j^{d+i} \psi_{t,1}(\boldsymbol{\theta}, \hat{\mathbf{h}}_t(\boldsymbol{\theta})) \partial_i^k \hat{\mathbf{h}}_t(\boldsymbol{\theta}), \qquad t \in \mathbb{N}, \end{aligned}$$

for indices $j \in \{1, \ldots, q\}, k \in \{1, \ldots, d\}$, or abridged, $\hat{\mathbf{h}}_{t+1}' = \hat{\tilde{\Phi}}_t(\hat{\mathbf{h}}_t')$. The replacement of $\hat{\mathbf{h}}_t$ by $\mathbf{h}_t$ and $\hat{\mathbf{h}}_t'$ by $\mathbf{d}_t$ in (6.3) leads to the following (linear) SRE $\mathbf{d}_{t+1} = \tilde{\Phi}_t(\mathbf{d}_t)$ on $\mathbb{C}(K, \mathbb{R}^{dq})$:

$$(6.4) \qquad \begin{aligned}{} [\mathbf{d}_{t+1}(\boldsymbol{\theta})]_\ell &= \partial_j^k \psi_{t,1}(\boldsymbol{\theta}, \mathbf{h}_t(\boldsymbol{\theta})) \\ &\quad + \sum_{i=1}^q \partial_j^{d+i} \psi_{t,1}(\boldsymbol{\theta}, \mathbf{h}_t(\boldsymbol{\theta})) [\mathbf{d}_t(\boldsymbol{\theta})]_{(i-1)d+k}, \qquad t \in \mathbb{Z}, \end{aligned}$$

where $\ell = (j-1)\,d + k \in \{1, \ldots, dq\}$. We now show the following:

(1) The SRE (6.4) has a unique stationary solution $(\mathbf{d}_t)$, which is ergodic. The random element $\mathbf{d}_t$ is $\mathcal{F}_{t-1}$-measurable for every $t$.

(2) We have that $\|\mathbf{\dot{h}}'_t - \mathbf{d}_t\|_K \xrightarrow{\text{e.a.s.}} 0$ as $t \to \infty$, that is, $(\mathbf{d}_t)_{t \in \mathbb{N}}$ is a stationary approximation of $(\mathbf{\dot{h}}'_t)_{t \in \mathbb{N}}$.

(3) The random functions $\mathbf{h}_t$ are a.s. continuously differentiable on $K$, and for each $t \in \mathbb{Z}$,

$$\mathbf{d}_t \equiv \mathbf{h}'_t.$$

Since we establish relation (1) via Theorem 2.8, we need to show that $\dot{\Phi}_t^{(r)}$ is a contraction on average for $r$ large enough. As is obvious from elementary calculus,

$$
(6.5) \quad \begin{aligned}
[\dot{\Phi}_t^{(r)}(\mathbf{d})(\boldsymbol{\theta})]_\ell &= \partial_j^k \psi_{t,r}(\boldsymbol{\theta}, \mathbf{h}_{t-r+1}(\boldsymbol{\theta})) \\
&\quad + \sum_{i=1}^q \partial_j^{d+i} \psi_{t,r}(\boldsymbol{\theta}, \mathbf{h}_{t-r+1}(\boldsymbol{\theta}))[\mathbf{d}(\boldsymbol{\theta})]_{(i-1)\,d+k}.
\end{aligned}
$$

From $\psi_{t,r}(\boldsymbol{\theta}, \mathbf{u}) = p_2 \circ \varphi_t^{(r)}(\boldsymbol{\theta}, \mathbf{u}) = \Phi_t^{(r)}(\mathbf{u})$, we deduce that

$$|\psi_{t,r}(\boldsymbol{\theta}, \mathbf{u}) - \psi_{t,r}(\boldsymbol{\theta}, \tilde{\mathbf{u}})| \le \Lambda(\Phi_t^{(r)})|\mathbf{u} - \tilde{\mathbf{u}}|, \qquad \mathbf{u}, \tilde{\mathbf{u}} \in [0, \infty)^q,$$

and therefore,

$$(6.6) \quad \sup_{\boldsymbol{\theta} \in K} |\partial_j^{d+i} \psi_{t,r}(\boldsymbol{\theta}, \mathbf{h}_{t-r+1}(\boldsymbol{\theta}))| \le \Lambda(\Phi_t^{(r)}), \qquad t \in \mathbb{Z},$$

for all $i \in \{1, \ldots, q\}$. Using the representation (6.5) and applying inequality (6.6), we obtain

$$
\begin{aligned}
&|[\dot{\Phi}_t^{(r)}(\mathbf{d})(\boldsymbol{\theta})]_\ell - [\dot{\Phi}_t^{(r)}(\tilde{\mathbf{d}})(\boldsymbol{\theta})]_\ell| \\
&\qquad \le \Lambda(\Phi_t^{(r)}) \sum_{i=1}^q |[\mathbf{d}(\boldsymbol{\theta})]_{(i-1)\,d+k} - [\tilde{\mathbf{d}}(\boldsymbol{\theta})]_{(i-1)\,d+k}| \\
&\qquad \le \text{const} \times \Lambda(\Phi_t^{(r)})\|\mathbf{d} - \tilde{\mathbf{d}}\|_K, \qquad \mathbf{d}, \tilde{\mathbf{d}} \in \mathbb{C}(K, \mathbb{R}^{dq}),
\end{aligned}
$$

for all $\ell \in \{1, \ldots, dq\}$ and $\boldsymbol{\theta} \in K$, whence

$$\Lambda(\dot{\Phi}_t^{(r)}) \le c\Lambda(\Phi_t^{(r)})$$

for a certain constant $c > 0$ not depending on $r$. Since $\mathbb{E}[\log \Lambda(\Phi_0^{(r)})] \to -\infty$ as $r \to \infty$, we can choose $r$ so large that

$$\mathbb{E}[\log \Lambda(\dot{\Phi}_0^{(r)})] \le \log c + \mathbb{E}[\log \Lambda(\Phi_0^{(r)})] < 0.$$

Thus, the SRE $\mathbf{d}_{t+1} = \dot{\Phi}_t(\mathbf{d}_t)$, $t \in \mathbb{Z}$, obeys condition S.2 of Theorem 2.8, and S.1 is true by virtue of (6.1). Consequently, the latter SRE admits a unique stationary ergodic solution $(\mathbf{d}_t)$, for which $\mathbf{d}_t$ is $\mathcal{F}_{t-1}$-measurable for every $t$. As regards the limit relation (2), we need to study the perturbed SRE

$$\mathbf{q}_{t+1} = \hat{\dot{\Phi}}_t(\mathbf{q}_t), \qquad t \in \mathbb{N},$$

which has $(\hat{\mathbf{h}}_t')_{t \in \mathbb{N}}$ as one of its solutions. By the assumption (6.2), the triangle inequality, Proposition 3.12 and an application of Lemma 2.1, we have that

$$\|\hat{\dot{\Phi}}_t(0) - \dot{\Phi}_t(0)\|_K \leq \text{const} \times \bar{C}_1(t) \|\hat{\mathbf{h}}_t - \mathbf{h}_t\|_K^\kappa \xrightarrow{\text{e.a.s.}} 0, \qquad t \to \infty,$$

and

$$\Lambda(\hat{\dot{\Phi}}_t - \dot{\Phi}_t) \leq \text{const} \times \bar{C}_1(t) \|\hat{\mathbf{h}}_t - \mathbf{h}_t\|_K^\kappa \xrightarrow{\text{e.a.s.}} 0, \qquad t \to \infty.$$

Now an application of Theorem 2.10 demonstrates $\|\hat{\mathbf{h}}_t' - \mathbf{d}_t\|_K \xrightarrow{\text{e.a.s.}} 0$ as $t \to \infty$. It remains to prove relation (3). Since $K$ coincides with the closure of its interior, the continuous differentiability of $\mathbf{h}_t$ on $K$ can be established by showing the existence of a sequence $(\mathbf{f}_n)_{n \in \mathbb{N}}$ of continuously differentiable functions on $K$ such that $\mathbf{f}_n \xrightarrow{\text{a.s.}} \mathbf{h}_t$ in $\mathbb{C}(K, \mathbb{R}^q)$ and $\mathbf{f}_n' \xrightarrow{\text{a.s.}} \mathbf{d}_t$ in $\mathbb{C}(K, \mathbb{R}^{dq})$ as $n \to \infty$; see Theorem 5.9.12 in [21]. For every *fixed* $m \geq 0$, the sequence $((\mathbf{h}_{t,m}', \mathbf{d}_t))$ is stationary ergodic by virtue of Proposition 2.5 [see (3.29) for the definition of $\mathbf{h}_{t,m}$]. Since $\mathbf{h}_{m,m}' = \hat{\mathbf{h}}_m'$, another application of Proposition 2.5 implies $\mathbf{h}_{t,m}' - \mathbf{d}_t \overset{\text{d}}{=} \hat{\mathbf{h}}_m' - \mathbf{d}_m$, $m \geq 0$. On the other hand, we have already shown that $\|\hat{\mathbf{h}}_m' - \mathbf{d}_m\|_K \xrightarrow{\text{e.a.s.}} 0$. Thus, $\|\mathbf{h}_{t,m}' - \mathbf{d}_t\|_K \xrightarrow{\mathbb{P}} 0$ as $m \to \infty$, and therefore, there is a subsequence $\mathbf{h}_{t,m_n}'$ with $\|\mathbf{h}_{t,m_n}' - \mathbf{d}_t\|_K \xrightarrow{\text{a.s.}} 0$ as $n \to \infty$. If we set $\mathbf{f}_n = \mathbf{h}_{t,m_n}$, $n \in \mathbb{N}$, then the sequence $(\mathbf{f}_n)_{n \in \mathbb{N}}$ satisfies $\mathbf{f}_n \xrightarrow{\text{a.s.}} \mathbf{h}_t$ in $\mathbb{C}(K, \mathbb{R}^q)$ and $\mathbf{f}_n' \xrightarrow{\text{a.s.}} \mathbf{d}_t$ in $\mathbb{C}(K, \mathbb{R}^{dq})$ as $n \to \infty$. This completes the proof of assertion (3). Summarizing, we have obtained the following proposition.

PROPOSITION 6.1. *Assume that conditions* D.1 *and* D.2 *are satisfied. Then the SRE* $\mathbf{d}_{t+1} = \dot{\Phi}_t(\mathbf{d}_t)$ *defined by* (6.4) *has a unique stationary solution* $(\mathbf{d}_t)$, *which is ergodic. For every* $t \in \mathbb{Z}$, *the random element* $\mathbf{d}_t$ *is* $\mathcal{F}_{t-1}$-*measurable. For every* $t \in \mathbb{Z}$, *the first derivatives of* $\mathbf{h}_t$ *coincide with* $\mathbf{s}_t$ *on* $K$ *a.s. Moreover,*

$$\|\hat{\mathbf{h}}_t' - \mathbf{d}_t\|_K \xrightarrow{\text{e.a.s.}} 0, \qquad t \to \infty.$$

*This justifies the following definition for the first derivatives of* $\mathbf{h}_t$: $\mathbf{h}_t' \equiv \mathbf{d}_t$.

Similar results can be derived for $(\hat{\mathbf{h}}_t'')_{t\in\mathbb{N}}$. At the origin of the analysis we have the observation that $(\hat{\mathbf{h}}_t'')_{t\in\mathbb{N}}$ obeys an SRE, which is contractive provided certain regularity assumptions hold. One can more or less follow the lines of proof of Proposition 6.1. In addition to D.1 and D.2, we will also assume the following set of conditions.

D.3 The function $(\boldsymbol{\theta}, \mathbf{s}) \mapsto g_{\boldsymbol{\theta}}(\mathbf{x}, \mathbf{s})$ on $K \times [0, \infty)^q$ is twice continuously differentiable for every fixed $\mathbf{x} \in \mathbb{R}^p$. For all $j \in \{1, \ldots, q\}$ and $k_1, k_2 \in \{1, \ldots, d+q\}$,

$$(6.7) \qquad \mathbb{E}\left[ \log^+\left( \sup_{\boldsymbol{\theta} \in K} |\partial_j^{k_1, k_2} \psi_{0,1}(\boldsymbol{\theta}, \mathbf{h}_0(\boldsymbol{\theta}))| \right) \right] < \infty.$$

The sequence of first derivatives $(\mathbf{h}_t')$ satisfies $\mathbb{E}(\log^+ \|\mathbf{h}_0'\|_K) < \infty$. Moreover, there exist a stationary sequence $(\bar{C}_2(t))$ with $\mathbb{E}[\log^+ \bar{C}_2(0)] < \infty$ and $\tilde{\kappa} \in (0, 1]$ such that

$$(6.8) \qquad \begin{aligned} \sup_{\boldsymbol{\theta} \in K} |\partial_j^{k_1, k_2} &\psi_{t,1}(\boldsymbol{\theta}, \mathbf{u}) - \partial_j^{k_1, k_2} \psi_{t,1}(\boldsymbol{\theta}, \tilde{\mathbf{u}})| \\ &\leq \bar{C}_2(t) |\mathbf{u} - \tilde{\mathbf{u}}|^{\tilde{\kappa}}, \qquad \mathbf{u}, \tilde{\mathbf{u}} \in [0, \infty)^q, \end{aligned}$$

for every $j \in \{1, \ldots, q\}$, $k_1, k_2 \in \{1, \ldots, d+q\}$ and $t \in \mathbb{Z}$.

PROPOSITION 6.2. *Assume that conditions D.1–D.3 are satisfied. Then*

$$\|\hat{\mathbf{h}}_t'' - \mathbf{e}_t\|_K \overset{\text{e.a.s.}}{\longrightarrow} 0, \qquad t \to \infty,$$

*where $(\mathbf{e}_t)$ is characterized as the unique stationary ergodic solution of the (linear) SRE (6.9) below. For every $t \in \mathbb{Z}$, the random element $\mathbf{e}_t$ is $\mathcal{F}_{t-1}$-measurable. The second derivatives of $\mathbf{h}_t$ on $K$ coincide with $\mathbf{e}_t$ a.s. for every $t$. Therefore, the following definition for the second derivatives of $\mathbf{h}_t$ is justified: $\mathbf{h}_t'' \equiv \mathbf{e}_t$.*

PROOF. Differentiation of both sides of (6.3) with respect to $\boldsymbol{\theta}$ shows that $\hat{\mathbf{h}}_{t+1}'' = \dot{\tilde{\Phi}}_t(\hat{\mathbf{h}}_t'')$, where $\dot{\tilde{\Phi}}_t$ is a linear random map. More precisely, for every $j \in \{1, \ldots, q\}$ and $k_1, k_2 \in \{1, \ldots, d\}$,

$$\partial_j^{k_1, k_2} \hat{\mathbf{h}}_{t+1}(\boldsymbol{\theta}) = \partial_j^{k_1, k_2} \psi_{t,1}(\boldsymbol{\theta}, \hat{\mathbf{h}}_t(\boldsymbol{\theta})) + \sum_{i=1}^q \partial_j^{k_2, d+i} \psi_{t,1}(\boldsymbol{\theta}, \hat{\mathbf{h}}_t(\boldsymbol{\theta}))\, \partial_i^{k_1} \hat{\mathbf{h}}_t(\boldsymbol{\theta})$$

$$+ \sum_{i=1}^q \Bigg( \partial_j^{k_1, d+i} \psi_{t,1}(\boldsymbol{\theta}, \hat{\mathbf{h}}_t(\boldsymbol{\theta}))$$

$$+ \sum_{i'=1}^q \partial_j^{d+i', d+i} \psi_{t,1}(\boldsymbol{\theta}, \hat{\mathbf{h}}_t(\boldsymbol{\theta})) \partial_{i'}^{k_1} \hat{\mathbf{h}}_t(\boldsymbol{\theta}) \Bigg) \partial_i^{k_2} \hat{\mathbf{h}}_t(\boldsymbol{\theta})$$

$$+ \sum_{i=1}^q \partial_j^{d+i} \psi_{t,1}(\boldsymbol{\theta}, \hat{\mathbf{h}}_t(\boldsymbol{\theta}))\, \partial_i^{k_1, k_2} \hat{\mathbf{h}}_t(\boldsymbol{\theta}), \qquad t \in \mathbb{N}.$$

This suggests to consider the (linear) SRE on $\mathbb{C}(K, \mathbb{R}^{d^2 q})$,

$$(6.9) \qquad \mathbf{e}_{t+1} = \ddot{\bar{\Phi}}_t(\mathbf{e}_t), \qquad t \in \mathbb{Z},$$

where for $\ell = (j-1)\,d^2 + (k_1 - 1)\,d + k_2 \in \{1, \ldots, d^2 q\}$,

$$[\ddot{\bar{\Phi}}_t(\mathbf{e})(\boldsymbol{\theta})]_\ell = \partial_j^{k_1, k_2} \psi_{t,1}(\boldsymbol{\theta}, \mathbf{h}_t(\boldsymbol{\theta})) + \sum_{i=1}^q \partial_j^{k_2, d+i} \psi_{t,1}(\boldsymbol{\theta}, \mathbf{h}_t(\boldsymbol{\theta}))\, \partial_i^{k_1} \mathbf{h}_t(\boldsymbol{\theta})$$

$$+ \sum_{i=1}^q \Bigg( \partial_j^{k_1, d+i} \psi_{t,1}(\boldsymbol{\theta}, \mathbf{h}_t(\boldsymbol{\theta}))$$

$$+ \sum_{i'=1}^q \partial_j^{d+i', d+i} \psi_{t,1}(\boldsymbol{\theta}, \mathbf{h}_t(\boldsymbol{\theta}))\, \partial_{i'}^{k_1} \mathbf{h}_t(\boldsymbol{\theta}) \Bigg) \partial_i^{k_2} \mathbf{h}_t(\boldsymbol{\theta})$$

$$+ \sum_{i=1}^q \partial_j^{d+i} \psi_{t,1}(\boldsymbol{\theta}, \mathbf{h}_t(\boldsymbol{\theta}))[\mathbf{e}(\boldsymbol{\theta})]_{(i-1)\,d^2 + (k_1 - 1)\,d + k_2},$$

$$\mathbf{e} \in \mathbb{C}(K, \mathbb{R}^{d^2 q}).$$

Exploiting the conditions (6.7) and $\mathbb{E}(\log^+ \|\mathbf{h}_0'\|_K) < \infty$ in D.3, one shows with exactly the same arguments as in the proof of Proposition 6.1 that the SRE (6.9) obeys the conditions of Theorem 2.8, which implies that it has a unique stationary ergodic solution $(\mathbf{e}_t)$. Furthermore, $\mathbf{e}_t$ is $\mathcal{F}_{t-1}$-measurable for every $t$. By means of the decomposition

$$\hat{a}\hat{b}\hat{c} - abc = (\hat{a} - a)bc + (\hat{a} - a)(\hat{b} - b)c + (\hat{a} - a)(\hat{b} - b)(\hat{c} - c)$$

$$+ (\hat{a} - a)b(\hat{c} - c) + a(\hat{b} - b)c + a(\hat{b} - b)(\hat{c} - c) + ab(\hat{c} - c)$$

and application of the bounds (6.8) together with $\|\hat{\mathbf{h}}_t - \mathbf{h}_t\|_K^{\tilde{\kappa}} \xrightarrow{\text{e.a.s.}} 0$ and $\|\hat{\mathbf{h}}_t' - \mathbf{h}_t'\|_K^{\tilde{\kappa}} \xrightarrow{\text{e.a.s.}} 0$, it can be verified that $(\ddot{\bar{\Phi}}_t)$ and $(\hat{\ddot{\bar{\Phi}}}_t)_{t \in \mathbb{N}}$ satisfy the conditions of Theorem 2.10. Thus, $\|\hat{\mathbf{h}}_t'' - \mathbf{e}_t\|_K \xrightarrow{\text{e.a.s.}} 0$. Analogously to the proof of Proposition 6.1, one demonstrates $\mathbf{e}_t \equiv \mathbf{h}_t''$. This concludes the proof. $\quad\square$

## 7. Asymptotic normality of the QMLE.
As we indicated in the Introduction, it is convenient to first establish the asymptotic normality of $\tilde{\boldsymbol{\theta}}_n$, defined in (4.4), and second to establish the asymptotic equivalence of $\tilde{\boldsymbol{\theta}}_n$ and $\hat{\boldsymbol{\theta}}_n$, that is, $\sqrt{n}(\hat{\boldsymbol{\theta}}_n - \tilde{\boldsymbol{\theta}}_n) \xrightarrow{\text{a.s.}} 0$. Following the classical approach, we will establish the asymptotic normality of $\tilde{\boldsymbol{\theta}}_n$ by means of a Taylor expansion of $L_n' = (\sum_{t=1}^n \ell_t)'$, where $\ell_t = -2^{-1}(\log h_t + X_t^2/h_t)$. For this reason, it is essential to study the limit properties of $L_n'$ and $L_n''$. Now we formulate the basic assumptions used throughout this section:

N.1 The assumptions C.1–C.4 of Section 4 are satisfied fulfilled and the true parameter $\boldsymbol{\theta}_0$ lies in the interior of the compact set $K$.

N.2 The assumptions D.1–D.3 of Proposition 6.2 are met so that $h_t$ is twice continuously differentiable on $K$.

N.3 The following moment conditions hold:

(i) $\mathbb{E} Z_0^4 < \infty$,

(ii) $\mathbb{E}\left(\dfrac{|h_0'(\boldsymbol{\theta}_0)|^2}{\sigma_0^4}\right) < \infty$,

(iii) $\mathbb{E}\|\ell_0'\|_K < \infty$,

(iv) $\mathbb{E}\|\ell_0''\|_K < \infty$.

N.4 The components of the vector $\frac{\partial g_\theta}{\partial \boldsymbol{\theta}}(\mathbf{X}_0, \boldsymbol{\sigma}_0^2)|_{\boldsymbol{\theta}=\boldsymbol{\theta}_0}$ are linearly independent random variables.

By virtue of Theorem 4.1, condition N.1 implies consistency. The requirement that $\boldsymbol{\theta}_0$ is in the interior of $K$ and N.2 enable the differentiation of $L_n$ in an open neighborhood of $\boldsymbol{\theta}_0$. Assumption N.4 assures that the asymptotic covariance matrix is regular. Observe that, with probability one,

$$(7.1) \qquad \ell_t'(\boldsymbol{\theta}) = -\frac{1}{2}\frac{h_t'(\boldsymbol{\theta})}{h_t(\boldsymbol{\theta})}\left(1 - \frac{X_t^2}{h_t(\boldsymbol{\theta})}\right),$$

$$(7.2) \qquad \ell_t''(\boldsymbol{\theta}) = -\frac{1}{2}\frac{1}{h_t(\boldsymbol{\theta})^2}\left((h_t'(\boldsymbol{\theta}))^T h_t'(\boldsymbol{\theta})\left(2\frac{X_t^2}{h_t(\boldsymbol{\theta})} - 1\right) + h_t''(\boldsymbol{\theta})(h_t(\boldsymbol{\theta}) - X_t^2)\right).$$

From Propositions 3.1, 3.12, 6.1 and 6.2, we infer that $(\ell_t')$ and $(\ell_t'')$ are stationary ergodic sequences of random elements with values in $\mathbb{C}(K, \mathbb{R}^d)$ and $\mathbb{C}(K, \mathbb{R}^{d \times d})$, respectively. An inspection of the proof of Theorem 4.1 shows that condition N.1 also implies $\tilde{\boldsymbol{\theta}}_n \xrightarrow{\text{a.s.}} \boldsymbol{\theta}_0$. Consequently, for large enough $n$, the following Taylor expansion is valid:

$$(7.3) \qquad L_n'(\tilde{\boldsymbol{\theta}}_n) = L_n'(\boldsymbol{\theta}_0) + L_n''(\boldsymbol{\zeta}_n)(\tilde{\boldsymbol{\theta}}_n - \boldsymbol{\theta}_0),$$

where $|\boldsymbol{\zeta}_n - \boldsymbol{\theta}_0| < |\tilde{\boldsymbol{\theta}}_n - \boldsymbol{\theta}_0|$. Since $\tilde{\boldsymbol{\theta}}_n$ is the maximizer of $L_n$ and $\boldsymbol{\theta}_0$ lies in the interior of $K$, one has $L_n'(\tilde{\boldsymbol{\theta}}_n) = 0$. Therefore, (7.3) is equivalent to

$$(7.4) \qquad n^{-1} L_n''(\boldsymbol{\zeta}_n)(\tilde{\boldsymbol{\theta}}_n - \boldsymbol{\theta}_0) = -n^{-1} L_n'(\boldsymbol{\theta}_0).$$

On account of $\mathbb{E}\|\ell_0''\|_K < \infty$ and the stationarity and ergodicity of $(\ell_t'')$, we may apply Theorem 2.7 to obtain $L_n''/n \xrightarrow{\text{a.s.}} L''$ in $\mathbb{C}(K, \mathbb{R}^{d \times d})$ as $n \to \infty$, where $L''(\boldsymbol{\theta}) = \mathbb{E}[\ell_0''(\boldsymbol{\theta})]$, $\boldsymbol{\theta} \in K$. This uniform convergence result together with $\boldsymbol{\zeta}_n \xrightarrow{\text{a.s.}} \boldsymbol{\theta}_0$ implies

$$L_n''(\boldsymbol{\zeta}_n)/n \xrightarrow{\text{a.s.}} \mathbb{E}[\ell_0''(\boldsymbol{\theta}_0)] = \mathbf{B}_0, \qquad n \to \infty.$$

By Propositions 3.12, 6.1 and 6.2, $h_0$, $h_0'$ and $h_0''$ are $\mathcal{F}_{-1}$-measurable. Exploiting $h_0(\boldsymbol{\theta}_0) = \sigma_0^2$ a.s., $X_0 = \sigma_0 Z_0$ and the independence of $Z_0$ and $\mathcal{F}_{-1}$, one may conclude that

$$(7.5) \qquad \mathbf{B}_0 = -2^{-1}\mathbb{E}[(h_0'(\boldsymbol{\theta}_0))^T h_0'(\boldsymbol{\theta}_0)/\sigma_0^4].$$

It is shown in Lemma 7.2 below that $\mathbf{B}_0$ is invertible. Consequently, the matrix $L_n''(\boldsymbol{\zeta}_n)/n$ has inverse $\mathbf{B}_0^{-1}(1 + o_{\mathbb{P}}(1))$, $n \to \infty$, and (7.4) has

$$\sqrt{n}(\tilde{\boldsymbol{\theta}}_n - \boldsymbol{\theta}_0) = -\mathbf{B}_0^{-1}(1 + o_{\mathbb{P}}(1))L_n'(\boldsymbol{\theta}_0)/\sqrt{n}, \qquad n \to \infty,$$

as a consequence. Therefore, the limit of $\sqrt{n}(\tilde{\boldsymbol{\theta}}_n - \boldsymbol{\theta}_0)$ is determined by that of $-\mathbf{B}_0^{-1} L_n'(\boldsymbol{\theta}_0)/\sqrt{n}$. Since $h_t(\boldsymbol{\theta}_0) = \sigma_t^2$ a.s. and $X_t = \sigma_t Z_t$,

$$L_n'(\boldsymbol{\theta}_0) = \sum_{t=1}^n \ell_t'(\boldsymbol{\theta}_0) = \frac{1}{2}\sum_{t=1}^n \frac{h_t'(\boldsymbol{\theta}_0)}{\sigma_t^2}(Z_t^2 - 1).$$

Since the random element $h_t'/\sigma_t^2$ is $\mathcal{F}_{t-1}$-measurable and since $\mathcal{F}_{t-1}$ is independent of $Z_t$ and $\mathbb{E}Z_t^2 = 1$, the sequence $(\ell_t'(\boldsymbol{\theta}_0))_{t\in\mathbb{N}}$ is a stationary ergodic zero-mean martingale difference sequence with respect to the filtration $(\mathcal{F}_t)_{t\in\mathbb{N}}$. By virtue of the moment condition N.3, the sequence $(\ell_t'(\boldsymbol{\theta}_0))_{t\in\mathbb{N}}$ is furthermore square integrable. Consequently, we can apply the central limit theorem for square-integrable stationary ergodic martingale difference sequences; see Theorem 18.3 in [3], which says that $n^{-1/2}L_n'(\boldsymbol{\theta}_0) \xrightarrow{\mathrm{d}} \mathcal{N}(\mathbf{0}, \mathbb{E}[(\ell_0'(\boldsymbol{\theta}_0))^T \ell_0'(\boldsymbol{\theta}_0)])$, $n \to \infty$. Together with (7.5), we conclude

$$\sqrt{n}(\tilde{\boldsymbol{\theta}}_n - \boldsymbol{\theta}_0) \xrightarrow{\mathrm{d}} \mathcal{N}(\mathbf{0}, \mathbf{V}_0), \qquad n \to \infty,$$

where

$$(7.6) \qquad \mathbf{V}_0 = 4^{-1}\mathbb{E}(Z_0^4 - 1)(\mathbb{E}[(h_0'(\boldsymbol{\theta}_0))^T h_0'(\boldsymbol{\theta}_0)/\sigma_0^4])^{-1}.$$

It is shown in Lemma 7.4 below that $\sqrt{n}|\hat{\boldsymbol{\theta}}_n - \tilde{\boldsymbol{\theta}}_n| \xrightarrow{\mathrm{a.s.}} 0$ so that an application of Slutsky's lemma finalizes the proof of the following theorem.

THEOREM 7.1. *Under the conditions* N.1–N.4, *the QMLE* $\hat{\boldsymbol{\theta}}_n$ *is strongly consistent and asymptotically normal, that is,*

$$\sqrt{n}(\hat{\boldsymbol{\theta}}_n - \boldsymbol{\theta}_0) \xrightarrow{\mathrm{d}} \mathcal{N}(\mathbf{0}, \mathbf{V}_0), \qquad n \to \infty,$$

*where the asymptotic covariance matrix* $\mathbf{V}_0$ *is given by* (7.6).

First we show that the asymptotic covariance matrix is regular.

LEMMA 7.2. *The assumptions* N.1–N.4 *imply that* $\mathbf{B}_0 = \mathbb{E}[\ell_0''(\boldsymbol{\theta}_0)]$ *is negative definite.*

PROOF. $\mathbf{B}_0$ being negative definite is equivalent to $\mathbf{C}_0 = \mathbb{E}[(h_0'(\boldsymbol{\theta}_0))^T h_0'(\boldsymbol{\theta}_0)/\sigma_0^4]$ being positive definite. It is evident that $\mathbf{C}_0$ is positive semi-definite. Assume $\mathbf{x}_0^T \mathbf{C}_0 \mathbf{x}_0 = 0$, for some $\mathbf{x}_0 \in \mathbb{R}^d$. This is equivalent to

$$\mathbb{E}\left|\frac{h_0'(\boldsymbol{\theta}_0)\mathbf{x}_0}{\sigma_0^2}\right|^2 = 0,$$

which implies $h_0'(\boldsymbol{\theta}_0)\mathbf{x}_0 = 0$ a.s., and also $\mathbf{h}_0'(\boldsymbol{\theta}_0)\mathbf{x}_0 = 0$ due to the stationarity of $(\mathbf{h}_t')$. Multiplying the equation

$$h_1'(\boldsymbol{\theta}_0) = \frac{\partial g_{\boldsymbol{\theta}}}{\partial \boldsymbol{\theta}}(\mathbf{X}_0, \boldsymbol{\sigma}_0^2)\Big|_{\boldsymbol{\theta}=\boldsymbol{\theta}_0} + \frac{\partial g_{\boldsymbol{\theta}}}{\partial \mathbf{s}}(\mathbf{X}_0, \boldsymbol{\sigma}_0^2)\Big|_{\boldsymbol{\theta}=\boldsymbol{\theta}_0} \mathbf{h}_0'(\boldsymbol{\theta}_0)$$

from the right by $\mathbf{x}_0$ and accounting for $h_1'(\boldsymbol{\theta}_0)\mathbf{x}_0 = 0$ results in $(\partial g_{\boldsymbol{\theta}}(\mathbf{X}_0, \boldsymbol{\sigma}_0^2)/\partial \boldsymbol{\theta})|_{\boldsymbol{\theta}=\boldsymbol{\theta}_0}\mathbf{x}_0 = 0$ a.s. Condition N.4 implies $\mathbf{x}_0 = 0$. This concludes the proof. $\square$

The remaining steps are devoted to the proof of Lemma 7.4.

LEMMA 7.3. *The assumptions* N.1–N.2 *imply*

$$n^{-1/2}\|\hat{L}_n' - L_n'\|_K \xrightarrow{\text{a.s.}} 0, \qquad n \to \infty.$$

PROOF. Notice that C.3 implies $\hat{h}_t(\boldsymbol{\theta}), h_t(\boldsymbol{\theta}) \geq \underline{g} > 0$ for all $\boldsymbol{\theta} \in K$. This and the mean value theorem applied to the function $f(a,b) = ab^{-1}(1 - X_t^2/b)$, $a \in \mathbb{R}$, $b \geq \underline{g}$, yield

$$\begin{aligned}
\|\hat{\ell}_t' - \ell_t'\|_K &= \left\|\frac{\hat{h}_t'}{\hat{h}_t}\left(1 - \frac{X_t^2}{\hat{h}_t}\right) - \frac{h_t'}{h_t}\left(1 - \frac{X_t^2}{h_t}\right)\right\|_K \\
&\leq C(1 + X_t^2)\{\|\hat{h}_t - h_t\|_K + \|\hat{h}_t' - h_t'\|_K + \|\hat{h}_t' - h_t'\|_K^2\|h_t'\|_K\}
\end{aligned}$$
(7.7)

for some $C > 0$. Recall that we assume $\mathbb{E}(\log^+ \sigma_0^2) < \infty$, $\mathbb{E}(\log^+ \|h_0'\|_K) < \infty$, and observe that Propositions 3.12 and 6.1 imply $\|\hat{h}_t - h_t\|_K \xrightarrow{\text{e.a.s.}} 0$ and $\|\hat{h}_t' - h_t'\|_K \xrightarrow{\text{e.a.s.}} 0$ as $t \to \infty$. Now (7.7) together with an application of Lemmas 2.1 and 2.2 shows $\|\hat{L}_n' - L_n'\|_K \leq \sum_{t=1}^{\infty}\|\hat{\ell}_t' - \ell_t'\|_K < \infty$ a.s. This completes the proof. $\square$

LEMMA 7.4. *The assumptions* N.1–N.4 *imply*

$$\sqrt{n}|\tilde{\boldsymbol{\theta}}_n - \hat{\boldsymbol{\theta}}_n| \xrightarrow{\text{a.s.}} 0, \qquad n \to \infty.$$

PROOF. From the mean value theorem,

$$L_n'(\tilde{\boldsymbol{\theta}}_n) - L_n'(\hat{\boldsymbol{\theta}}_n) = L_n''(\tilde{\boldsymbol{\zeta}}_n)(\tilde{\boldsymbol{\theta}}_n - \hat{\boldsymbol{\theta}}_n),$$
(7.8)

where $\tilde{\boldsymbol{\zeta}}_n$ lies on the line segment connecting $\hat{\boldsymbol{\theta}}_n$ and $\tilde{\boldsymbol{\theta}}_n$. This line segment is completely contained in the interior of $K$ provided $n$ is large enough. Since $L_n'(\hat{\boldsymbol{\theta}}_n) = \hat{L}_n'(\hat{\boldsymbol{\theta}}_n) = 0$, equation (7.8) is equivalent to

$$n^{-1/2}(\hat{L}_n'(\hat{\boldsymbol{\theta}}_n) - L_n'(\hat{\boldsymbol{\theta}}_n)) = n^{-1}L_n''(\tilde{\boldsymbol{\zeta}}_n)n^{1/2}(\tilde{\boldsymbol{\theta}}_n - \hat{\boldsymbol{\theta}}_n).$$
(7.9)

By virtue of Lemma 7.3, both sides of (7.9) tend to 0 a.s. as $n \to \infty$. Because of N.3(iv), we can apply Theorem 2.7 to $L_n''/n$ and together with $\tilde{\boldsymbol{\zeta}}_n \xrightarrow{\text{a.s.}} \boldsymbol{\theta}_0$ conclude $L''_n(\tilde{\boldsymbol{\zeta}}_n)/n \xrightarrow{\text{a.s.}} \mathbf{B}_0 = \mathbb{E}[\ell''_0(\boldsymbol{\theta}_0)]$. Since the matrix $\mathbf{B}_0$ is invertible, as shown by Lemma 7.2, we can deduce $\sqrt{n}(\hat{\boldsymbol{\theta}}_n - \tilde{\boldsymbol{\theta}}_n) \xrightarrow{\text{a.s.}} 0$, which completes the proof. $\square$

REMARK 7.5.   In general, it seems impossible to find a tractable expression for the asymptotic covariance matrix $\mathbf{V}_0$ due to the fact that the joint distribution of $(\sigma_0^2, h'_0(\boldsymbol{\theta}_0))$ is not known, not even for GARCH(1, 1). It is, however, possible to consistently estimate $\mathbf{V}_0$ from the data. Defining the residual by $\hat{Z}_t^{(n)} = X_t/(\hat{h}_t(\hat{\boldsymbol{\theta}}_n))^{1/2}$,

$$\hat{\mathbf{V}}_0^{(n)} = \left( \frac{1}{4n} \sum_{t=1}^n (\hat{Z}_t^4 - 1) \right) \left( \frac{1}{n} \sum_{t=1}^n \frac{(\hat{h}'_t(\hat{\boldsymbol{\theta}}_n))^T \hat{h}'_t(\hat{\boldsymbol{\theta}}_n)}{\hat{h}_t(\hat{\boldsymbol{\theta}}_n)^2} \right)^{-1}$$

is a strongly consistent estimator for the matrix $\mathbf{V}_0$. We sketch how this can be demonstrated. Define

$$\mathbf{M}_n(\boldsymbol{\theta}) = \frac{1}{n} \sum_{t=1}^n \frac{(h'_t(\boldsymbol{\theta}))^T h'_t(\boldsymbol{\theta})}{(h_t(\boldsymbol{\theta}))^2} \quad \text{and}$$

$$\hat{\mathbf{M}}_n(\boldsymbol{\theta}) = \frac{1}{n} \sum_{t=1}^n \frac{(\hat{h}'_t(\boldsymbol{\theta}))^T \hat{h}'_t(\boldsymbol{\theta})}{(\hat{h}_t(\boldsymbol{\theta}))^2}, \qquad \boldsymbol{\theta} \in K,$$

and suppose $\mathbb{E}\|h'_0/h_0\|_K^2 < \infty$ in addition to N.1–N.4. Then by an application of Theorem 2.7, $\mathbf{M}_n \xrightarrow{\text{a.s.}} \mathbf{M}$ in $\mathbb{C}(K, \mathbb{R}^{d \times d})$, where $\mathbf{M}(\boldsymbol{\theta}) = \mathbb{E}[(h'_0(\boldsymbol{\theta}))^T h'_0(\boldsymbol{\theta})/ (h_0(\boldsymbol{\theta}))^2]$. Using the same method as in Lemma 7.3, one derives a bound for $\|\hat{\mathbf{M}}_n - \mathbf{M}_n\|_K$ and shows $\|\hat{\mathbf{M}}_n - \mathbf{M}_n\|_K \xrightarrow{\text{a.s.}} 0$, which implies $\hat{\mathbf{M}}_n \xrightarrow{\text{a.s.}} \mathbf{M}$ in $\mathbb{C}(K, \mathbb{R}^{d \times d})$ as $n \to \infty$. Therefore, $(\hat{\mathbf{M}}_n(\hat{\boldsymbol{\theta}}_n))^{-1} \xrightarrow{\text{a.s.}} (\mathbf{M}(\boldsymbol{\theta}_0))^{-1}$. Likewise, $n^{-1} \sum_{t=1}^n ((\hat{Z}_t^{(n)})^4 - 1)$ can be shown to converge to $\mathbb{E}(Z_0^4 - 1)$ a.s. We also mention that the practical implementation for the computation of the matrix $\hat{\mathbf{V}}_0^{(n)}$ becomes particularly simple if one makes use of the recursion for $\hat{\mathbf{h}}'_t$; see equation (6.3).

## 8. Examples: asymptotic normality of QMLE in AGARCH($p, q$).
At the moment we cannot provide a proof of the asymptotic normality of the QMLE in the general EGARCH model, because it is intricate to verify the moment conditions in N.3. However, one can deal with the special case when $\beta_0 = 0$; see [31]. In this paper we concentrate on establishing the asymptotic normality of the QMLE in AGARCH($p, q$) and thereby generalize the results of [2].

Take a compact set $K \subset (0, \infty) \times [0, \infty)^p \times B \times [-1, 1]$ which coincides with the closure of its interior. Assume that the true parameter vector $\boldsymbol{\theta}_0 = (\alpha_0^\circ, \alpha_1^\circ, \ldots, \alpha_p^\circ, \beta_1^\circ, \ldots, \beta_q^\circ, \gamma^\circ)^T$ is contained in the interior of $K$ and suppose the conditions of Theorem 5.5 hold. Analogously to [2], we assume $\mathbb{E}Z_0^4 < \infty$ and suppose there is a $\mu > 0$ such that

$$(8.1) \qquad\qquad \mathbb{P}(|Z_0| \le z) = o(z^\mu), \qquad z \downarrow 0.$$

The verification of conditions D.1–D.3 of Section 6 is rather straightforward. The only steps which require some care are the moment conditions. The arguments in Remark 3.6 imply that there is a $\tilde{q} > 0$ with $\mathbb{E}|\sigma_0|^{2\tilde{q}} < \infty$, and hence, also $\mathbb{E}|X_0|^{2\tilde{q}} < \infty$. The Minkowski inequality applied to the left-hand side of (5.6) shows that $\mathbb{E}\|h_0\|_K^{\tilde{q}} < \infty$. Altogether, $\mathbb{E}(\log^+ \sigma_0^2) < \infty$ and $\mathbb{E}(\log^+ \|h_0\|_K) < \infty$. By standard methods, one can also show that the differential and the series in (5.6) can be interchanged. Then $\mathbb{E}(\|h_0'\|_K^{\tilde{q}}) < \infty$ and $\mathbb{E}(\log^+ \|h_0'\|_K^{\tilde{q}}) < \infty$ follow. Altogether we have established N.1 and N.2. It is less easy to prove the moment conditions of N.3. At the origin lies the observation that the random variables $\|h_0'/h_0\|_K$ and $\|h_0''/h_0\|_K$ have finite moments of any order and that

$$\mathbb{E}\|X_0^2/h_0\|_K^\nu < \infty$$

for any $\nu < 2$; this follows from the same ideas and techniques as used to prove Lemmas 5.1 and 5.2 in [2]. Therefore, by an application of the Hölder inequality to the norm of (7.1), $\mathbb{E}\|\ell_0'\|_K < \infty$. As to the second derivative $\ell_0''$, we mention that the inequality

$$\|xy^T\|_K \le \|x\|_K \|y\|_K$$

for any two elements $x, y \in \mathbb{C}(K, \mathbb{R}^d)$ together with the triangle inequality applied to (7.2) implies

$$\|\ell_0''\|_K \le \frac{1}{2}\left(\left\|\frac{h_0'}{h_0}\right\|_K^2 \left(2\left\|\frac{X_0^2}{h_0}\right\|_K + 1\right) + \left\|\frac{h_0''}{h_0}\right\|_K \left(1 + \left\|\frac{X_0^2}{h_0}\right\|_K\right)\right).$$

By an application of the Hölder inequality, $\mathbb{E}\|\ell_0''\|_K < \infty$. Thus, conditions N.1–N.3 of Section 7 are satisfied. N.4 will be verified in Lemma 8.2 below. Now an application of Theorem 7.1 yields the following result.

THEOREM 8.1.  *Let $(X_t)$ be a stationary AGARCH$(p,q)$ process with true parameter vector $\boldsymbol{\theta}_0 = (\alpha_0^\circ, \alpha_1^\circ, \ldots, \alpha_p^\circ, \beta_1^\circ, \ldots, \beta_q^\circ, \gamma^\circ)^T$ in the interior of the compact set $K \subset (0,\infty) \times [0,\infty)^p \times B \times [-1,1]$, where $B = \{(\beta_1, \ldots, \beta_q)^T \in [0,1)^q \mid \sum_{j=1}^q \beta_j < 1\}$. Suppose also that $K$ coincides with the closure of its interior. Assume that the conditions of Theorem 5.5 hold and that $Z_0$ satisfies $\mathbb{E}Z_0^4 < \infty$ and (8.1). Then the QMLE (4.2) is strongly consistent and asymptotically normal with asymptotic covariance matrix $\mathbf{V}_0$ given by (7.6).*

LEMMA 8.2.  *Under the conditions imposed by Theorem 8.1, the condition N.4 is satisfied.*

PROOF.    By straightforward computation,

$$\frac{\partial g_{\boldsymbol{\theta}}}{\partial \boldsymbol{\theta}}(\mathbf{X}_0, \boldsymbol{\sigma}_0^2)\Big|_{\boldsymbol{\theta}=\boldsymbol{\theta}_0}$$
$$= (1, (|X_0| - \gamma^\circ X_0)^2, \ldots, (|X_{-p+1}| - \gamma^\circ X_{-p+1})^2, \sigma_0^2, \ldots, \sigma_{-q+1}^2,$$
$$-2\alpha_1^\circ X_0(|X_0| - \gamma^\circ X_0) - \cdots - 2\alpha_p^\circ X_{-p+1}(|X_{-p+1}| - \gamma^\circ X_{-p+1}))^T.$$

Assume for a $\boldsymbol{\xi} = (\lambda_0, \ldots, \lambda_p, \mu_1, \ldots, \mu_q, \eta)^T \in \mathbb{R}^{p+q+2}$ that

$$(8.2) \qquad\qquad \frac{\partial g_{\boldsymbol{\theta}}}{\partial \boldsymbol{\theta}}(\mathbf{X}_0, \boldsymbol{\sigma}_0^2)\Big|_{\boldsymbol{\theta}=\boldsymbol{\theta}_0} \boldsymbol{\xi} = 0.$$

Writing out the latter equation results in

$$\lambda_0 + \sigma_0^2(\lambda_1(|Z_0| - \gamma^\circ Z_0)^2 + \mu_1 - 2\eta\alpha_1^\circ Z_0(|Z_0| - \gamma^\circ Z_0)) + Y_{-1} = 0 \qquad \text{a.s.},$$

where $Y_{-1}$ is a certain $\mathcal{F}_{-1}$-measurable random variable. Consequently, the multiplier of $\sigma_0^2$ in the latter identity must be degenerate, that is,

$$(8.3) \quad (\lambda_1(1 + (\gamma^\circ)^2) + 2\eta\alpha_1^\circ\gamma^\circ)Z_0^2 - (2\gamma^\circ\lambda_1 + 2\eta\alpha_1^\circ)Z_0|Z_0| + \mu_1 = c \qquad \text{a.s.}$$

for a certain $c \in \mathbb{R}$. Since the distribution of $Z_0$ is not concentrated at two points, $\{1, Z_0|Z_0|, Z_0^2\}$ are linearly independent. Combining this information with (8.3) and taking into account that $\alpha_1^\circ > 0$ and $|\gamma^\circ| < 1$, we conclude $\lambda_1 = \eta = \mu_1 - c = 0$. Thus, we must show the linear independence of $1, (|X_0| - \gamma^\circ X_0)^2, \ldots, (|X_{-p+1}| - \gamma^\circ X_{-p+1})^2, \sigma_0^2, \ldots, \sigma_{-q+1}^2$. The following arguments are similar to the ones used by Berkes, Horváth and Kokoszka [2] for the proof of their Lemma 5.7. Equation (8.2) with $\eta = 0$ and the stationarity of $((X_t, \sigma_t))$ imply

$$\lambda_0 + \sum_{i=1}^p \lambda_i(|X_{t-i}| - \gamma^\circ X_{t-i})^2 + \sum_{j=1}^q \mu_j\sigma_{t-j}^2 = 0 \qquad \text{a.s.},$$

which equivalently written in backshift operator notation becomes

$$(8.4) \qquad \lambda_0 + \lambda(B)(|X_t| - \gamma^\circ X_t)^2 + \mu(B)\sigma_t^2 = 0 \qquad \text{a.s.},$$

with

$$\lambda(z) = \sum_{i=1}^p \lambda_i z^i \quad \text{and} \quad \mu(z) = \sum_{j=1}^q \mu_j z^j.$$

Recall Lemma 5.2, where we derived the a.s. representation

$$(8.5) \qquad \sigma_t^2 = \frac{\alpha_0^\circ}{b^\circ(1)} + a^\circ(B)(b^\circ(B))^{-1}(|X_t| - \gamma^\circ X_t)^2$$

with $a^\circ(z) = \sum_{i=1}^p \alpha_i^\circ z^i$ and $b^\circ(z) = 1 - \sum_{j=1}^q \beta_j^\circ z^j$. Plugging (8.5) into (8.4), we obtain that

$$\lambda_0 + \frac{\mu(1)\alpha_0^\circ}{b^\circ(1)} + (\lambda(B) + \mu(B)a^\circ(B)(b^\circ(B))^{-1})(|X_t| - \gamma^\circ X_t)^2 = 0 \qquad \text{a.s.,}$$

and with the identical arguments as given in Lemma 5.4, we conclude

$$\lambda(z) + \mu(z)\frac{a^\circ(z)}{b^\circ(z)} \equiv 0 \quad \text{and} \quad \lambda_0 + \frac{\mu(1)\alpha_0^\circ}{b^\circ(1)} \equiv 0.$$

If $\mu \neq 0$, the rational function $\mu(z)a^\circ(z)/b^\circ(z)$ has at least $p+1$ zeros (counted with their multiplicities), whereas $\lambda(z)$ has at most $p$ zeros, which is a contradiction. Therefore, $\mu(z) = \lambda(z) = 0$ and $\boldsymbol{\xi} = 0$, which completes the proof.  □

## 9. Nonstationarities.
So far our considerations were based on the fundamental assumption that we sample from the unique *stationary* solution to the SRE (3.1) of the general heteroscedastic model. This is a strong assumption, particularly in light of Remark 3.2. There we mentioned that simulating the stationary solution is, in general, impossible, and we provided an algorithm which generates a sequence $({}_k\tilde{\boldsymbol{\sigma}}_t^2)_{t\in\mathbb{N}}$ which approaches the stationary solution $({}_k\boldsymbol{\sigma}_t^2)_{t\in\mathbb{N}}$ of (3.2) with an error decaying to zero exponentially fast as $t \to \infty$:

(1) Take an initial value $\boldsymbol{\varsigma}_0^2 \in [0, \infty)^k$, set ${}_k\tilde{\boldsymbol{\sigma}}_0^2 = \boldsymbol{\varsigma}_0^2$, and generate ${}_k\tilde{\boldsymbol{\sigma}}_t^2$ according to (3.2).

(2) Set $(\tilde{\sigma}_t^2, \tilde{X}_t) = ({}_k\tilde{\boldsymbol{\sigma}}_{t,1}^2)^{1/2}(1, Z_t)$, $t = 0, 1, \ldots$.

Notice that $\boldsymbol{\varsigma}_0 \in [0, \infty)^q$ is an arbitrarily chosen initial value. Correspondingly, we define $\tilde{X}_t = \tilde{\sigma}_t Z_t$ and observe that $|\tilde{X}_t^2 - X_t^2| = Z_t^2|\tilde{\sigma}_t^2 - \sigma_t^2| \overset{\text{e.a.s.}}{\longrightarrow} 0$ as $t \to \infty$ by virtue of Proposition 3.1 and Lemma 2.1. To the best of our knowledge, simulation studies have been based on random samples of $(\tilde{X}_t, t = -p+1, \ldots, n)$, rather than of $(X_t, t = -p+1, \ldots n)$; see, for example, [23]. The effects of this inherent nonstationarity due to some initialization error are hardly ever addressed. Although simulation studies are most often carried out for determining the *small sample* behavior of estimators, from a theoretical point of view, a minimal requirement to validate the simulation approach would be the asymptotic equivalence of the QMLE applied to $(\tilde{X}_t)_{t \geq -p+1}$ and to $(X_t)_{t \geq -p+1}$. Let $\tilde{L}_n = \sum_{t=1}^n (X_t^2/\tilde{h}_t + \log \tilde{h}_t)$ denote the likelihood based on the observations $(\tilde{X}_t, t = -p+1, \ldots, n)$. To establish the asymptotic equivalence of $\hat{\boldsymbol{\theta}}_n$ and of the maximizer of $\tilde{L}_n$, we merely need to show $n^{-1}\|\tilde{L}_n - \hat{L}_n\|_K \overset{\text{e.a.s.}}{\longrightarrow} 0$ and $n^{-1/2}\|\tilde{L}_n' - \hat{L}_n'\|_K \overset{\text{e.a.s.}}{\longrightarrow} 0$. One can impose conditions on the functions $g_{\boldsymbol{\theta}}$, which imply the latter two relations, but since everything is in line with the ideas already presented in this paper, we omit the details.

**Acknowledgments.**  The authors would like to thank an anonymous referee who called our attention to Theorem 3.1 of [5]. This enabled us to generalize the first version of this paper and to treat conditionally heteroscedastic processes of arbitrary order $(p, q)$. Furthermore, we thank Anders Rahbek (Copenhagen) and Friedrich Liese (Rostock) for helpful discussions. The paper contains material from the Ph.D. thesis of Daniel Straumann, which was written at the Laboratory of Actuarial Mathematics, Department of Applied Mathematics and Statistics, University of Copenhagen. The thesis appeared in the Springer *Lecture Notes in Statistics* Series [31]. Thomas Mikosch's research is partially supported by MaPhySto, The Danish Research Network for Mathematical Physics and Stochastics, and DYNSTOCH, a research training network under the program Improving Human Potential financed by the Fifth Framework Programme of the European Commission.

RISKMETRICS GROUP
GENEVA BUSINESS CENTER
12 AVENUE DES MORGINES
CH-1213 PETIT-LANCY
SWITZERLAND
E-MAIL: daniel.straumann@riskmetrics.com

LABORATORY OF ACTUARIAL MATHEMATICS
DEPARTMENT OF APPLIED MATHEMATICS
    AND STATISTICS
UNIVERSITY OF COPENHAGEN
UNIVERSITETSPARKEN 5
DK-2100 COPENHAGEN Ø
DENMARK
AND
MAPHYSTO
THE DANISH RESEARCH NETWORK
    FOR MATHEMATICAL PHYSICS
    AND STOCHASTICS
E-MAIL: mikosch@math.ku.dk